\documentclass[11pt]{article}
\usepackage[utf8]{inputenc}
\usepackage[
  a4paper,            
  textwidth=170mm,
  textheight=244.8mm,
  left=20mm,          
  right=20mm,         
  top=25mm,           
  headheight=11pt,
  headsep=22pt,
  footskip=37pt
]{geometry}

\usepackage{graphicx}
\usepackage{amsmath}
\usepackage{amssymb}
\usepackage[T1]{fontenc}
\usepackage{amsfonts}
\usepackage{graphicx}
\usepackage{lipsum}

\usepackage{algorithmic}
\usepackage{algorithm}
\usepackage{ulem}

\usepackage{hyperref}
\usepackage[capitalise]{cleveref}
\usepackage{nicematrix}

\Crefname{equation}{Equation}{Equations} 
\Crefname{figure}{Figure}{Figures}
\Crefname{section}{Section}{Sections}
\Crefname{subsection}{Sub-Section}{Sub-Sections}

\numberwithin{equation}{section}
\numberwithin{table}{section}
\numberwithin{figure}{section}
\numberwithin{algorithm}{section}

\newcommand{\pder}[2]{\frac{\partial#1}{\partial#2}}

\newcommand{\norm}[1]{\left\lVert#1\right\rVert}

\newcommand{\vare}[1]{{#1}_{e}}
\newcommand{\vars}[1]{{#1}_{s}}
\newcommand{\varcc}[1]{{#1}_{s}}
\newcommand{\parame}[1]{{#1}_{e}}
\newcommand{\params}[1]{{#1}_{s}}
\newcommand{\paramam}[1]{{#1}_{am}}
\newcommand{\paramcc}[1]{{#1}_{cc}}

\newcommand{\ndvare}[1]{\tilde{#1}_{e}}
\newcommand{\ndvars}[1]{\tilde{#1}_{s}}

\newcommand{\ndparamam}[1]{\tilde{#1}_{am}}
\newcommand{\ndparamcc}[1]{\tilde{#1}_{cc}}

\newcommand{\varei}[2]{{#1}_{e,\rm{#2}}}
\newcommand{\varsi}[2]{{#1}_{s,\rm{#2}}}

\newcommand{\cs}[1]{{#1}^{\star}}
\newcommand{\cse}[1]{{#1}^{\star}_{e}}
\newcommand{\css}[1]{{#1}^{\star}_{s}}

\newcommand{\paran}[1]{\left(#1\right)}
\newcommand{\curlybrac}[1]{\{#1\}}
\newcommand{\bigbrac}[1]{\left[#1\right]}

\newcommand{\auxvar}[3]{ {#1}_{#2}^{#3} }

\newcommand{\BV}[0]{Butler-Volmer }

\newcommand{\half}[0]{\frac{1}{2}}

\def\xs{\bar{x}}

\usepackage{authblk}

\title{High-order adaptive multi-domain time integration scheme for microscale lithium-ion batteries simulations}

\author[1,2]{Ali Asad \thanks{ali.asad@polytechnique.edu}}
\affil[1]{CMAP, CNRS, École polytechnique, Institut Polytechnique de Paris, 91120 Palaiseau Cedex, France}
\author[2]{Romain {de} Loubens}
\affil[2]{TotalEnergies OneTech, 3 Boulevard Thomas Gobert, 91120 Palaiseau Cedex, France}
\author[3]{Laurent François}
\affil[3]{ONERA, DMPE, 6 Chemin de la Vauve aux Granges, 91120 Palaiseau Cedex, France}
\author[1]{Marc Massot}

\begin{document}

\maketitle

\begin{abstract}
    We investigate the modeling and simulation of ionic transport and charge conservation in lithium-ion batteries (LIBs) at the microscale. It is a multiphysics problem that involves a wide range of time scales. 
    The associated computational challenges motivate the investigation of numerical techniques that can decouple the time integration of the governing equations in the liquid electrolyte and the solid phase (active materials and current collectors).
    First, it is shown that semi-discretization in space of the non-dimensionalized governing equations leads to a system of \textit{index-}1 \textit{semi-explicit} differential algebraic equations (DAEs).
    Then, a new generation of strategies for multi-domain integration is presented, enabling high-order adaptive coupling of both domains in time, with efficient and potentially different domain integrators. They reach a high level of flexibility for real applications, beyond the limitations of multirate methods.   
    A simple 1D LIB half-cell code is implemented as a demonstrator of the new strategy for the simulation of different modes of cell operation. 
    The integration of the decoupled subsystems is performed with high-order accurate implicit nonlinear solvers.
    The accuracy of the space discretization is assessed by comparing the numerical results to the analytical solutions.
    Then, temporal convergence studies demonstrate the accuracy of the new multi-domain coupling approach. 
    Finally, the accuracy and computational efficiency of the adaptive coupling strategy are discussed in the light of the conditioning of the decoupled subproblems compared to the one of the fully-coupled problem.
    This new approach will constitute a key ingredient for the high-fidelity 3D LIB simulations based on actual electrode microstructures.     
\end{abstract}
    
\subsubsection*{Keywords}
Lithium-ion batteries, High-order time integration methods, Adaptive multi-domain integration scheme
    
\section{Introduction} \label{sec:introduction}
The shift in energy paradigm is evident and, as the world moves towards more renewable energy sources, the role of chemical energy storage using batteries becomes more vital.
The rechargeable lithium-ion batteries (LIBs) with high energy and power densities are one of the widely used energy storage devices \cite{xue_lithium-ion_2016}.
They are especially used in electric transportation and aerospace applications due to their portability, robustness, and reliability. 
In recent years, the demand has also been growing for stationary storage applications, notably in the areas of renewable energy production and grid frequency regulation.
Given the established popularity of LIBs and the increase in demand, the improvement in LIB design and performance is an appealing and growing area of research. In this context, numerical simulation is an important tool, besides experimental studies, to gain insight into the functioning of LIBs. 
The reliability of such simulations depends on our ability to model LIBs mathematically and to devise accurate and efficient numerical methods. 

The LIB operation is governed by physical phenomena at different spatial and temporal scales. The mathematical modeling of LIB yields a multiphysics and multiscale problem. From an engineering perspective, the macroscale of interest is defined by the size of a battery cell or a stack of cells. Macroscale models, such as volume-averaged \cite{fuller_simulation_1994, doyle_use_1995} or equivalent circuit models \cite{liaw2004modeling}, benefit from their simplicity and are widely used in practice, e.g., for the electric and thermal control of LIBs \cite{di2010lithium,bartlett2015electrochemical}. However, they suffer from severe limitations, such as loss of resolution due to volume averaging and failure to capture the effect of heterogeneities or defects in the electrode microstructure \cite{latz_multiscale_2015, muller2018quantifying}. Better resolution can be achieved with microscale models, accounting for the real microstructure of porous electrodes \cite{kashkooli_multiscale_2016, hutzenlaub_three-dimensional_2014}, or with nanoscale models that can represent complex interfacial phenomena, like Solid Electrolyte Interphase (SEI) formation in graphite anodes \cite{muralidharan_molecular_2018, single_dynamics_2016}. Over the last two decades, the development of such descriptions has benefited from the improvement of high-resolution imaging and high-performance computing. As a result, microscale models now offer a good compromise between complexity and accuracy at the electrode scale.

Microscale LIB models are based on a micro-continuum description of transport phenomena and electrochemistry. In these models, each point belongs to a specific material domain, such as electrolyte, active material, separator, or current collector. Within each domain, governing equations may be written for the conservation of mass, electric charge, and energy. At each material interface, specific conditions should be supplied, such as the Butler-Volmer equation modeling the charge transfer at the solid-electrolyte interface \cite{latz_thermodynamic_2011}. In the liquid electrolyte domain, the mathematical description of ionic transport is given by the Nernst-Planck equation for dilute solutions or the Maxwell-Stefan theory for concentrated solutions \cite{psaltis_comparing_2011}. In general, the expression of charge conservation assumes solution electroneutrality, except near the electric double layer (EDL). However, in the scope of microscale models, EDL effects can be neglected as the pore size is usually much larger than the Debye length. In each active material domain, the transport of lithium is commonly represented by Fick's law and electronic conduction by Ohm's law. Such description relies on an effective diffusion coefficient that can be measured experimentally under representative operating conditions. It is ideally suited for active materials with solid-solution behavior, but also applicable with a lesser degree of accuracy to phase-separating materials like graphite and lithium iron phosphate (LFP).

In their standard formulation, microscale LIB models rely on the previously mentioned assumptions, while also neglecting thermal, mechanical and degradation effects. Their domain of validity is thus dependent on the nature of LIB materials and the operating conditions. In a number of cases, however, a standard microscale formulation is sufficient to reproduce experimental charge-discharge curves within acceptable accuracy \cite{less_micro-scale_2012, allen_segregated_2021,zhang_comparison_2014}. For more complex situations, extended formulations are available, including for example thermal effects \cite{latz_thermodynamic_2011, latz_multiscale_2015}, active material deformations \cite{bower_finite_2011, rahani_role_2013, ferraro_electrode_2020}, phase transformations \cite{castelli2019numerical}, lithium plating \cite{hein2020electrochemical} or SEI growth \cite{schneider2022efficient}.

The numerical solution of microscale LIB models can be obtained with various simulation techniques already adopted in other fields of computational physics. In \cite{popov_finite_2011,spotnitz2012geometry}, a discretization of the LIB microscale equations based on the finite volume method is presented. It is implemented in BEST\textit{micro} \cite{best_software}, a commercial LIB simulator using Cartesian meshes, as well as in the CFD-based simulation software STAR-CCM$+$\textsuperscript{\textregistered} using unstructured meshes \cite{ccm0u}. A simulation package for the microscale study of LIBs is also available in the finite-element software COMSOL Multiphysics\textsuperscript{\textregistered} \cite{comsol_software}. Besides commercial tools, several research simulators have been developed to model the microscale physics of LIBs based on the finite-element method \cite{roberts_framework_2014,Smith_2009}.

Simulation of LIB cells at the microscale presents multiple computational challenges. The spatial discretization of the governing equations yields a stiff system of equations, which is commonly seen in multiphysics and multiscale problems. Moreover, the \BV current interface condition introduces a strong nonlinearity into the system. 
Last but not least, 3D simulation of LIBs based on real electrode microstructures requires  large number of degrees of freedom (DOFs) and inversion of large ill-conditioned linear systems. 

An efficient way to solve numerically stiff systems of nonlinear equations is to use stable implicit solvers~\cite{hairer_solving_1996, ivp_dae_book}. However, such solvers are complex to implement and require costly Newton iterations, all the more because of the ill-conditioning of the related linear systems to be solved. 
While multirate methods initiated some time ago \cite{gear1980} within the framework of ordinary differential equations, analyzed and improved further in \cite{kvaerno2001}, can handle various time scales not only in applications where several operators are involved \cite{woodward2024} but also in coupled problems \cite{constantinescu2023}, they still require a monolithic time integration through an algebraic resolution \cite{constantinescu2022} of the resulting systems and involve constraints in the time integration of the various subdomains through a global approach. We aim at keeping a high level of flexibility in the way we handle the various domains in terms of both discretization and methods (or even models \cite{francois_multiphysical_2022}) and focus on methods where the various domains are handled separately.
In~\cite{iliev_domain_2016, iliev_splitting_2017}, three algorithms of time integration with domain splitting are proposed to reduce the number of Newton iterations. Convergence and stability of these algorithms are also studied for a 2D LIB model problem. In~\cite{GOLDIN2012118}, another splitting approach is implemented, where lithium concentrations and electric potentials are solved independently and subsequently coupled using Picard outer-iterations. 
A mortar-based spatial coupling scheme allows the use of different finite-element meshes in the electrolyte and electrode phases in~\cite{fang_monolthic_2018}, while using a monolithic time integration approach.  
In~\cite{allen_segregated_2021}, a domain-splitting technique in time is complemented with the block Gau\ss-Seidel (BGS) method and algebraic multigrid (AMG) for matrix inversion to solve the 3D LIB problem. The performance of these methods is studied in comparison with a direct linear solver.
Linear solvers relying on block preconditioners corresponding to the concentration and potential fields show improved computational efficiency~\cite{fang_parallel_2019}.
Overcoming the numerous challenges associated with the microscale simulation of LIBs remains an active area of research. Most of the present studies neglect the physical complexities of LIBs, such as active material phase transformations, lithium plating, mechanical effects, etc. 
Before adding more physics to the LIB model, there is a strong need for innovative strategies to  conduct high-fidelity simulations with a high level of accuracy and robustness, but also at a reasonable computational cost. 

In this work, we introduce and assess a novel multi-domain technique based on high-order adaptive coupling in time. The objective is to cope with the strong multiscale character of the LIB model as efficiently as possible, while maintaining a coupling strategy with excellent properties of stability and accuracy. The novelty relies on the introduction of a coupling technique, which is adaptive and high-order in time. The purpose of our contribution is threefold. First, we make a precise link between a representative half-cell LIB modeling and the mathematical structure of the semi-discretized system of partial differential equations in one dimension, which includes most of the difficulties, we will have to cope with in the multi-dimensional system. We formulate the problem under the form of a system of \textit{index-}1 \textit{semi-explicit} differential algebraic equations (DAEs), for which accurate and stable time integrators exist.
Second, we present the novel high-order and adaptive multi-domain time-coupling strategy and its analysis.
Third, the numerical strategy is assessed after a series of verification cases provided in \cref{sec:validation_result}. 
We show that the linear systems involved in the segregated domain problems are better conditioned than those of the monolithic fully-coupled approach. This indicates a strong potential of our method for 3D LIB cell high-fidelity simulations based on real electrode microstructures.

The paper is organized as follows. We introduce the 1D mathematical model for LIB half-cells in \cref{sec:lib_intro} and analyze the structure of the semi-discretized in space system of equations. We describe the numerical methods used for solving the LIB problem, including the multi-domain high-order adaptive coupling in time, in \cref{sec:num_schemes}. Furthermore, we present and discuss our results in \cref{sec:results_discussions}, which assess the strategy. Finally, relying on a study of the conditioning provided in \cref{sec:cnstudy_subsys},
we conclude in \cref{sec:conclusions} on the potential of the approach for realistic 3D simulations.

\section{1D LIB half-cell model} \label{sec:lib_intro}
In this section, we present the governing equations of the 1D LIB half-cell model along with its initial and boundary conditions.
Introducing characteristic scales, we derive a non-dimensional formulation. We then apply a spatial discretization on a finite volume grid. Finally, we describe the mathematical nature of the discretized system of equations.

\subsection{The microscale continuum description} \label{subsec:govn_equations}
\begin{figure}[htbp]
  \centering
  \includegraphics[width=0.9\linewidth]{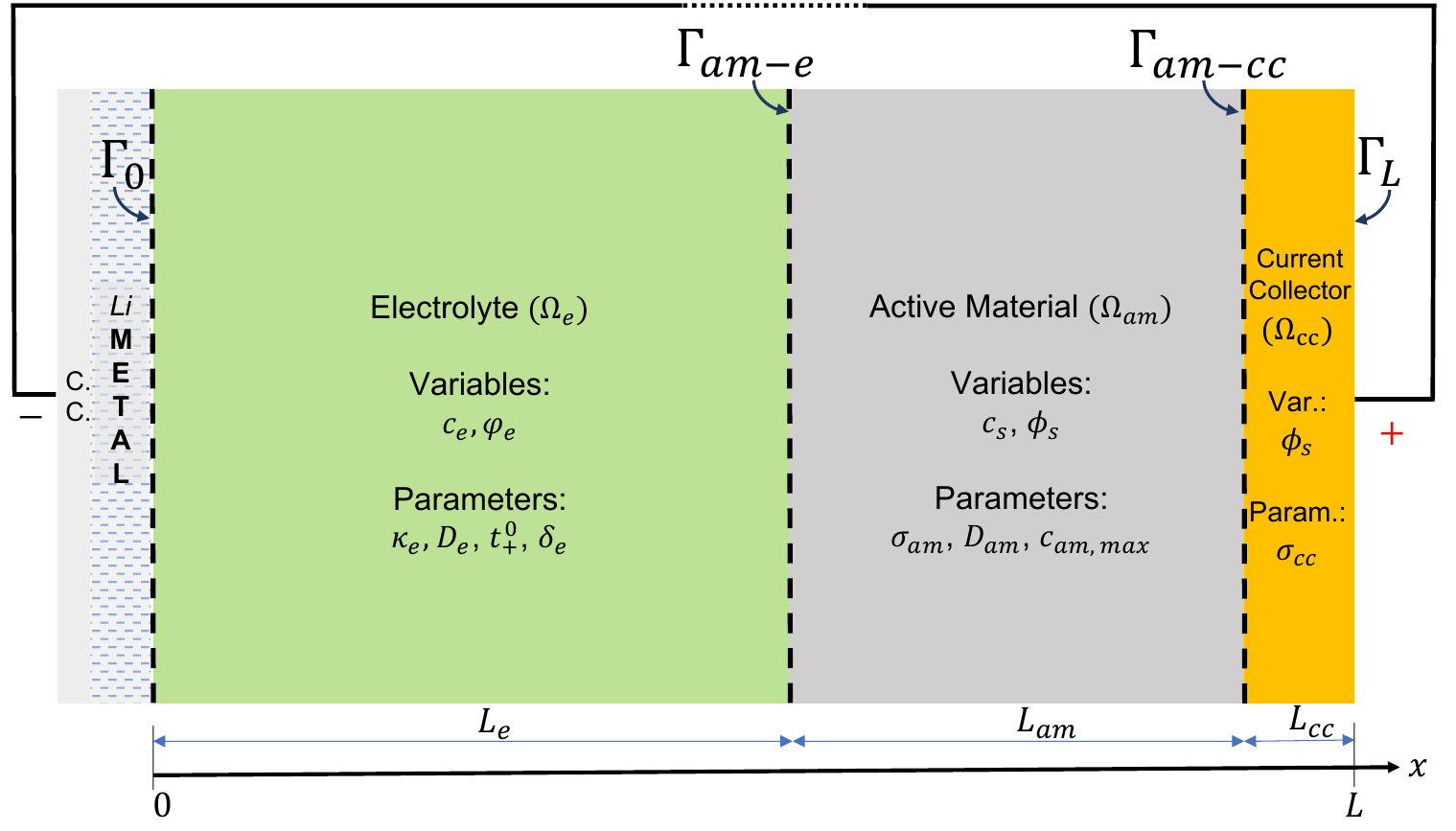}
  \caption{Schematic diagram of the 1D LIB half-cell model.}
  \label{fig:1d-schematic-halfcell}
\end{figure}

The LIB operation is governed by coupled equations of mass and charge conservation in the domains of electrolyte, active material, and current collector. A schematic diagram of a 1D half-cell is shown in \cref{fig:1d-schematic-halfcell}. Here the anode consists of lithium metal, which can be simply modeled as a boundary condition on the electrolyte variables (see~\cref{Eq:BC-Ne-at-interface-0} below). Compared to a full-cell model, the half-cell model is expected to exhibit similar numerical complexity. Indeed, the anode and cathode of the full-cell model are governed by the same equations (although with different parameters). We will now look at the governing equations and parameters in each domain as well as the boundary and initial conditions. 

\subsubsection{Electrolyte} \label{ssubsec:electrolyte}
The electrolyte domain $\Omega_{e}$ represents the liquid electrolyte of the LIB cell, assumed here to consist
of a single lithium salt. Using the electroneutrality assumption and the concentrated solution theory, the following equations for mass and charge conservation can be derived from the Stefan-Maxwell equations \cite{latz_thermodynamic_2011,richardson2022charge,bizeray2016resolving,newman2021electrochemical}    defining lithium-ion molar flux $\vare{\mathbf{N}}$ and ionic current density $\vare{\mathbf{i}}$ 
\begin{align}
	\pder{\vare{c}}{t} &= -\nabla\cdot\vare{\mathbf{N}} , \quad \vare{\mathbf{N}}\paran{\vare{c}, \vare{\varphi}} = - \parame{D} \nabla \vare{c} +  \frac{t_+^0}{F} \vare{\mathbf{i}},  \label{Eq:Electrolyte-Mass-Conservation} \\
	0 &= -\nabla\cdot\vare{\mathbf{i}} , \quad \vare{\mathbf{i}}\paran{\vare{c}, \vare{\varphi}} = - \parame{\kappa}\nabla\vare{\varphi} + \frac{2 RT(1-t_+^0)}{F}(1+\parame{\delta})\parame{\kappa}\nabla \ln \vare{c}. \label{Eq:Electrolyte-Charge-Conservation}
\end{align}
Here $F$ is Faraday's constant, $R$ is the ideal gas constant, and $T$ is the ambient temperature.   
The independent variables $\vare{c}$ and $\vare{\varphi}$ represent lithium-ion concentration and electrochemical potential, respectively. The electrochemical potential is defined as $\vare{\varphi} = \vare{\mu}(c_e)/F + \vare{\Phi}$, which is composed of electric potential $\vare{\Phi}$ and chemical potential $\vare{\mu}$ of the electrolyte.  
Model parameters include the transference number $t_+^0$, the electrolyte diffusivity coefficient $\parame{D}$, the ionic conductivity $\parame{\kappa}$, and the logarithmic derivative of the activity coefficient with respect to concentration, $\parame{\delta}$. This coefficient appears due to the classical relation between chemical potential and ionic concentration. In this work, the concentration dependency of these parameters is neglected.

\subsubsection{Electrode} \label{ssubsec:electrode}
The electrode (solid) domain $\params{\Omega}$ consists of the active material domain $\paramam{\Omega}$ and the current collector domain $\paramcc{\Omega}$. 
The governing equations for lithium transport and electric charge conservation in $\paramam{\Omega}$, as well as 
lithium molar flux $\vars{\mathbf{N}}$ and electric current density $\vars{\mathbf{i}}$, are given by

\begin{align}
	\pder{\vars{c}}{t} &= -\nabla\cdot\vars{\mathbf{N}}, \quad \vars{\mathbf{N}}\paran{\vars{c}} = - \paramam{D}\nabla{\vars{c}},\label{Eq:Solid-Lithium-Conservation} \\
	0 &= -\nabla\cdot\vars{\mathbf{i}}, \quad \vars{\mathbf{i}}\paran{\vars{\phi}} = - \paramam{\sigma}\nabla\vars{\phi}.\label{Eq:Solid-Charge-Conservation}
\end{align}
Here independent variables are lithium concentration $\vars{c}$ and electric potential $\vars{\phi}$. Model parameters include the active material diffusivity coefficient $\paramam{D}$ and electric conductivity $\paramam{\sigma}$, which are considered to be constant. We note that, similar to the approach in \cite{latz_thermodynamic_2011,richardson2012multiscale}, lithium is assumed to be in a neutral state, and the electric current is purely electronic. Alternative modeling approaches have been explored to avoid such assumptions \cite{farkhondeh2011mathematical,zeng2014phase}, particularly for phase-change electrode materials.

The current collector does not allow lithium intercalation or de-intercalation. Therefore, it is simply described by the charge conservation of~\cref{Eq:Solid-Charge-Conservation}, where the electric current density $\vars{\mathbf{i}}$ is now given by
\begin{align}
	\vars{\mathbf{i}}\paran{\vars{\phi}} = -\paramcc{\sigma}\nabla\vars{\phi} . \label{Eq:CurrentCollector-Current-Flux}
\end{align}
Here $\paramcc{\sigma}$ is the electric conductivity of the current collector. Typically, $\paramcc{\sigma}\gg\paramam{\sigma}$, i.e., the current collector is a much better electronic conductor than the active material. 

\subsubsection{Boundary and initial conditions} \label{bc-and-ic}
In \cref{fig:1d-schematic-halfcell}, external boundaries of the cell located at $x=0$ and $x=L$ are denoted as $\Gamma_{0}$ and $\Gamma_L$, respectively. Further, internal boundaries exist at interfaces between material domains, between electrolyte and active material ($\Gamma_{am-e}$), and between active material and current collector ($\Gamma_{am-cc}$). The boundary conditions at these interfaces are often referred to as interface conditions.    

The lithium exchange (intercalation/de-intercalation) between electrolyte and electrodes takes place at $\Gamma_{0}$ and $\Gamma_{am-e}$ via the redox reaction $\rm{Li}^+ + \rm{e}^- \rightleftharpoons \rm{Li}$. This process is modeled as an effective exchange current density at the interface which is described by the \BV model \cite{batchelor-mcauley_recent_2015}. It is given by
\begin{align}
    i_{BV} &= i_{0} \paran{e^{\frac{\alpha F}{RT}\eta} - e^{-\frac{(1-\alpha)F}{RT}\eta}} , \label{i-BV}
\end{align}
where $i_0$ denotes the equilibrium exchange current density, $\alpha$ is the anodic strength coefficient and $\eta$ is the electrode overpotential defined by 
\begin{align}
    \eta &= \vars{\phi}-\vare{\varphi}-U_0(\vars{c}) . \label{eta}
\end{align}
Here $U_0$ is the open-circuit potential of the electrode, which is a function of lithium concentration in active material. At $\Gamma_{am-e}$, $i_0 = F k_0 \sqrt{\vare{c} \vars{c}\paran{c_{s, \rm{max}}-\vars{c}}}$, where $k_0$ is the reaction rate constant and $c_{s, \rm{max}}$ is the maximum lithium concentration in active material. 
At $\Gamma_0$, the lithium metal anode leads to a simplified interface condition, with $\vars{\phi}=$~0 and  $U_0=$~0, hence $\eta=-\vare{\varphi}$. Moreover, the equilibrium exchange current density of lithium metal is assumed to be constant. It is denoted as $i_{0,\rm{Li}}$. 

In our study, we consider symmetric interface redox reactions, i.e. $\alpha =$~0.5. The resulting current densities at the anode ($\Gamma_0$) and cathode ($\Gamma_{am-e}$) interfaces are
\begin{align}
    i_{BV}^{(A)} &=  2 i_{0,\rm{Li}} \sinh \paran{-\frac{F}{2RT}\vare{\varphi}} , \label{i-BV-A} \\
    i_{BV}^{(C)} &=  2Fk_0 \sqrt{\vare{c} \vars{c}\paran{c_{s, \rm{max}}-\vars{c}}} \sinh \paran{\frac{F}{2RT}\paran{\vars{\phi}-\vare{\varphi}-U_0(\vars{c})}} , \label{i-BV-C} 
\end{align}
where superscripts $(A)$ and $(C)$ denote the anode and cathode, respectively. 
Therefore, at $\Gamma_0$ 
 \begin{align}
 	\vare{N}\paran{\vare{c}, \vare{\varphi}} |_{\Gamma_0} &= \frac{i_{BV}^{(A)}}{F}, \quad  \vare{i}\paran{\vare{c}, \vare{\varphi}} |_{\Gamma_0} = i_{BV}^{(A)}, \label{Eq:BC-Ne-at-interface-0}
\end{align}
while at $\Gamma_{am-e}$
\begin{align}
	\vare{N}\paran{\vare{c}, \vare{\varphi}} |_{\Gamma_{am-e}} &= - \frac{i_{BV}^{(C)}}{F}, \quad  \vare{i}\paran{\vare{c}, \vare{\varphi}} |_{\Gamma_{am-e}} = - i_{BV}^{(C)}, \label{Eq:BC-Ne-at-interface-am-e} \\
    \vars{N}\paran{\vars{c}}|_{\Gamma_{am-e}} &= - \frac{i_{BV}^{(C)}}{F} , \quad
    \vars{i}\paran{\vars{\phi}}|_{\Gamma_{am-e}} = - i_{BV}^{(C)}. \label{Eq:BC-is-at-interface-am-e}
\end{align}
Note that the fluxes in \crefrange{Eq:BC-Ne-at-interface-0}{Eq:BC-is-at-interface-am-e} are scalars resulting from their dot products with 1D normal vectors.  

At $\Gamma_{am-cc}$, lithium flux is not allowed,  it is also assumed that electric potential and current density are continuous
\begin{align}
    \vars{N}\paran{\vars{c}} |_{\Gamma_{am-cc}} &= 0, \quad
    \vars{\phi}|_{\Gamma_{am-cc}^{(-)}} = \vars{\phi}|_{\Gamma_{am-cc}^{(+)}}, \quad 
    \vars{i}|_{\Gamma_{am-cc}^{(-)}} = \vars{i}|_{\Gamma_{am-cc}^{(+)}}, \label{Eq:BC-icc-iam-at-interface-cc/am}
\end{align}
where $(-)$ and $(+)$ signs represent the left and right neighborhoods of the interface, respectively. The above conditions mean that any contact resistance or capacitance effect between active material and current collector is neglected.

The external boundary condition at $\Gamma_L$ is defined by the mode of operation of the LIB cell. A constant current  (\textit{galvanostatic} or CC) charge$/$discharge condition imposes a constant current density at $\Gamma_L$.
It is usually defined as the product of a reference current $i_{1C}$ (``1C'' rate) by a \textit{C-rate} factor $\xi$. The definition of $i_{1C}$ is based on the battery charge capacity, i.e., it is the current required for a full charge of the battery in one hour, starting from a fully discharged state. 
Hence, the imposed current density in CC mode is  
\begin{align}
	\varcc{i}\paran{\vars{\phi}}|_{\Gamma_{L}} &= -\xi i_{1C} . \label{Eq:BC-icc-at-interface-L-CC}
\end{align}
The negative sign comes from the sign convention for electric current. Note that a positive $\xi$ indicates charging while a negative value is for discharging. 
Further, a constant voltage (\textit{potentiostatic} or CV) mode of operation means a fixed potential $\phi_{s, \rm{ext}}$ is imposed at $\Gamma_{L}$, i.e.
\begin{align}
    \vars{\phi}|_{\Gamma_{L}} &= \phi_{s, \rm{ext}} . \label{Eq:BC-phic-at-interface-L-CC}
\end{align}
In some cases, the imposed external potential can be a function of time. Such a situation will be considered in our numerical example of \cref{subsec:res_adaptive_coupling}.  

Finally, the mathematical model of LIB half-cells is closed by introducing the initial conditions. Hence, we have
\begin{align}
    \vare{c}(x,0) &= c_{e,\rm{I}} , \quad 
    \vars{c}(x,0) = c_{s,\rm{I}} , \label{Eq:IC-cs}   
\end{align} 
with $c_{e,\rm{I}}$ is the initial electrolyte concentration, and $c_{s,\rm{I}}$ the initial lithium concentration in active material (given in \cref{tab:constants_parameters}).

\subsection{Non-dimensionalization} \label{subsec:non-dim_govn_equations}
In the following, all characteristic scales of the problem are denoted by a star (${}^{\star}$) superscript.  
First, we define $\cs{L}=L$ as a characteristic length scale. As a potential scale, we use the ambient thermal voltage $\cs{\Phi}=RT/F$ in both the electrolyte and solid domains. 
The lithium-ion concentration in the electrolyte domain is scaled by its initial value, i.e., $\cse{c}=c_{e,\rm{I}}$, 
while the scale for lithium concentration in active material is chosen as $\css{c}=c_{s,\rm{max}}$. 
As a common timescale for the electrolyte and solid domains, we use the electrolyte diffusion timescale, hence $\cs{t} = (\cs{L})^2/D_e$. These scales are summarized in \Cref{tab:scales}, which also provides definitions for the flux scales $\cse{N}$, $\cse{i}$, $\css{N}$, and $\css{i}$.  

Using the characteristic scales defined above, we now rewrite the governing equations in terms of dimensionless variables, denoted by a tilde ($\tilde{\ }$) superscript. Therefore, in the electrolyte domain
\begin{align}
    \partial_{\tilde{t}}\ndvare{c} &= - \partial_{\tilde{x}}\ndvare{N}, \quad   \ndvare{N}\paran{\ndvare{c}, \ndvare{\varphi}} = -\partial_{\tilde{x}}\ndvare{c} + \frac{\cse{i}}{\cse{N}} \frac{t_+^0}{F} \ndvare{i}, \label{Eq:Electrolyte-Mass-Flux-nd} \\
    0 &= - \partial_{\tilde{x}}\ndvare{i}, \quad   \ndvare{i}\paran{\ndvare{c}, \ndvare{\varphi}} = -\partial_{\tilde{x}}\ndvare{\varphi} + \frac{\kappa_D}{\ndvare{c}}\partial_{\tilde{x}}\ndvare{c}. \label{Eq:Electrolyte-Current-Flux-nd}
\end{align}
Here $\kappa_D = 2(1-t_+^0)(1+\delta_e)$ is a dimensionless parameter characterizing the diffusive contribution to ionic current. The other dimensionless quantity $(\cse{i}t_+^0)/(\cse{N}F)$ represents the ratio of electric migration flux to diffusion flux. It is usually referred to as the electric Péclet number~\cite{Arunachalam_2015}.    
Similarly, in the solid domain, if we note $\varepsilon_D = D_{am} / D_e$ and $\varepsilon_{\sigma} = \sigma_{am} / \sigma_{cc}$, we have
\begin{align}
    \partial_{\tilde{t}}\ndvars{c} &= - \varepsilon_{D}\partial_{\tilde{x}}\ndvars{N}, \quad \ndvars{N}\paran{\ndvars{c}} =
    \begin{cases}
        - \partial_{\tilde{x}}\ndvars{c}  , & \quad \tilde{x} \in \ndparamam{\Omega},  \\
        \quad 0 , & \quad \tilde{x} \in \ndparamcc{\Omega}, 
    \end{cases} \label{Eq:Solid-Mass-Flux-nd} \\
    0 &= - \partial_{\tilde{x}}\ndvars{i}, \quad  \ndvars{i}(\ndvars{\phi}) = 
    \begin{cases}
        - \partial_{\tilde{x}}\ndvars{\phi}  , & \quad \tilde{x} \in \ndparamam{\Omega}, \\
        - \frac{1}{\varepsilon_{\sigma}} \partial_{\tilde{x}}\ndvars{\phi} , & \quad \tilde{x} \in \ndparamcc{\Omega}.
    \end{cases} \label{Eq:Solid-Current-Flux-nd} 
\end{align}

We first define a unit scale for the \BV current density: $\cs{i}_{BV} = \mathrm{1\ A/m^2}$. Hence, the boundary conditions at $\tilde{\Gamma}_{0}$ become $\ndvare{N}=\tilde{i}_{BV}^{(A)}/F\cse{N}$ and $\ndvare{i}=\tilde{i}_{BV}^{(A)}/\cse{i}$. Similarly, at $\tilde{\Gamma}_{am-e}$ the dimensionless interface conditions are $\ndvare{N}=-\tilde{i}_{BV}^{(C)}/F\cse{N}$, $\ndvare{i}=-\tilde{i}_{BV}^{(C)}/\cse{i}$, $\ndvars{N}=-\tilde{i}_{BV}^{(C)}/F\css{N}$, and $\ndvars{i}=-\tilde{i}_{BV}^{(C)}/\css{i}$. Further, $\ndvars{N}=0$ at  $\tilde{\Gamma}_{am-cc}$. The current density or the potential at $\tilde{\Gamma}_{L}$, depending on the LIB operating mode, is non-dimensionalized using the flux scale $\css{i}$ or the potential scale $\cs{\Phi}$, respectively. Finally, the initial conditions become $\ndvare{c}(\tilde{x},0)=$~1 and $\ndvars{c}(\tilde{x},0)=c_{s,\rm{I}}/\css{c}$.

The non-dimensionalization of the governing equations shows that the time derivatives of  concentrations in electrolyte and active material differ in magnitude by a factor $\varepsilon_D$ (see \cref{Eq:Solid-Mass-Flux-nd}). Using typical values of $D_{am}$ and $D_e$, we infer that $\varepsilon_D \in [\mathrm{10^{-5},10^{-3}}]$.
It should be noted that the coexistence of two timescales in the LIB microscale model contributes to the ill-conditioning of the problem. Hence, it motivates us to explore multi-domain techniques for time integration, where improved conditioning is expected for the subproblems in each domain due to decoupling. For simplicity, we henceforth omit the ``~$\widetilde{\cdot}$~'' notation for the non-dimensional variables.

\subsection{Semi-discretization in space} 
\label{subsec:semi-discretized}
\begin{figure}[tbhp]
  \centering
  \includegraphics[width=0.65\linewidth]{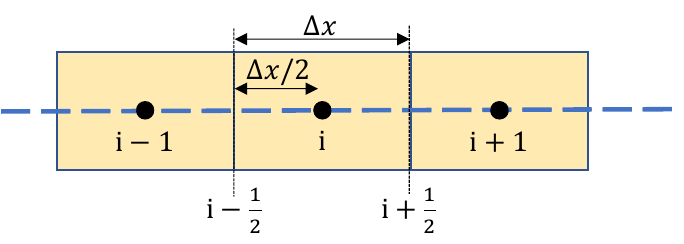}
  \caption{Schematic view of a 1D finite volume grid.}
  \label{fig:fv_cells}
\end{figure}

The finite volume method is well-suited for problems involving conservation laws such as heat and mass transfer \cite{EYMARD2000713}. In this work, the LIB half-cell model was derived from mass and charge conservation laws in each domain, with flux-based interface conditions. Thus, we perform the spatial discretization of the non-dimensional governing equations on a 1D finite volume grid, as illustrated in \cref{fig:fv_cells}. The discrete values of concentration and potential variables are located at cell centers, whereas molar fluxes and current densities are defined at cell faces. The former are denoted by subscript ${\rm i}$ and the latter by subscript ${\rm i+1/2}$, where ${\rm i}$ is the index of the ${\rm i}$-th cell.
To approximate the spatial derivatives of concentration and potential, we use the second-order accurate central differencing scheme. 

The discretized domains of electrolyte, active material, and current collector have $\parame{\mathcal{N}}$, $\paramam{\mathcal{N}}$, and $\paramcc{\mathcal{N}}$ finite volume cells, respectively. Hence, the number of cells in $\params{\Omega}$ is $\params{\mathcal{N}}=\paramam{\mathcal{N}}+\paramcc{\mathcal{N}}$, and the total number of cells $\mathcal{N}=\parame{\mathcal{N}}+\params{\mathcal{N}}$. Within a given material domain $\Omega_{\alpha}$ ($\alpha=e,am,cc$), we employ a constant cell width $\Delta x_{\alpha} = L_{\alpha}/\mathcal{N}_{\alpha}$ where $L_{\alpha}$ denotes the domain length. This ensures that interfaces between domains coincide with particular faces of the finite volume grid. 
In the following, we assume for simplicity but without loss of generality that the domain lengths allow to define a uniform grid size $\Delta x=\Delta x_{e}=\Delta x_{am}=\Delta x_{cc}$.

The system of equations obtained after spatial discretization in the electrolyte domain is given by
\begin{align}
    d_t \varei{c}{i} = \varei{F}{i} , \quad 0 = \varei{G}{i}, \ \label{Eqn:fv_electrolyte}
\end{align}
where, for $\text{i} = \mathrm{1,2,\ldots, \parame{\mathcal{N}}}$, we have 
\begin{align}
    \varei{F}{i} &= -\frac{1}{\Delta x}\paran{\varei{N}{i+\half} - \varei{N}{i-\half}}, \quad 
    \varei{G}{i} = -\frac{1}{\Delta x}\paran{\varei{i}{i+\half} - \varei{i}{i-\half}}. \label{Eqn:fv_electrolyte_current_flux}
\end{align}
Here, we obtain $\varei{N}{i+\half}$ and $\varei{i}{i+\half}$ from spatial discretization of the fluxes defined in \cref{Eq:Electrolyte-Mass-Flux-nd,Eq:Electrolyte-Current-Flux-nd}, respectively. We get
\begin{align}
    \varei{N}{i \pm \half} &= - \paran{ \pm\frac{\varei{c}{i \pm 1} - \varei{c}{i}}{\Delta x} }  + \frac{\cse{i}}{\cse{N}} \frac{t_+^0}{F} \varei{i}{\rm i \pm \half}, \label{Eq:fv-Electrolyte-Mass-Flux-nd}
\end{align}
and  
\begin{align}
    \varei{i}{ i \pm \half} &= - \paran{ \pm \frac{ \varei{\varphi}{i \pm 1} - \varei{\varphi}{i}}{\Delta x}} + \frac{\kappa_{D}}{\varei{c}{i \pm \half}} \paran{\pm \frac{ \varei{c}{i \pm 1} - \varei{c}{i}}{\Delta x}}, \label{Eq:fv-Electrolyte-Current-Flux-nd}
\end{align} 
for ${\rm i} = \mathrm{2,3,\ldots, \parame{\mathcal{N}}-1}$. The associated boundary conditions define the fluxes at ${\rm i} = 1$ and ${\rm i} = \parame{\mathcal{N}}$. To ensure accurate flux computation in \cref{Eq:fv-Electrolyte-Current-Flux-nd}, concentrations at cell faces ($\varei{c}{i \pm \half}$) are approximated using the harmonic mean of the adjacent cell-centered values.

Similarly, in the solid domain, we have
\begin{align}
    d_t \varsi{c}{j} = \varsi{F}{j}, \quad \quad 0 = \varsi{G}{j}, \label{Eqn:fv_solid}
\end{align}
where, for $\text{j} = \mathrm{1,2,\ldots, \params{\mathcal{N}}}$, we have
\begin{align}
    \varsi{F}{j} =& -\frac{\varepsilon_D}{\Delta x}\paran{\varsi{N}{j+\half} - \varsi{N}{j-\half}}, \quad
    \varsi{G}{j} = -\frac{1}{\Delta x}\paran{\varsi{i}{j+\half} - \varsi{i}{j-\half}}. \label{Eqn:fv_solid_current_flux}
\end{align}
The numerical fluxes $\varsi{N}{j+\half}$ and $\varsi{i}{j+\half}$ are the discrete versions of \cref{Eq:Solid-Mass-Flux-nd,Eq:Solid-Current-Flux-nd} given by  
\begin{align}
    \varsi{N}{i \pm \half} &= - \paran{ \pm\frac{\varsi{c}{i \pm 1} - \varsi{c}{i}}{\Delta x} }, \label{Eq:fv-Solid-Mass-Flux-nd}
\end{align}
and  
\begin{align}
    \varsi{i}{ i \pm \half} &=  - \varsi{\sigma}{i \pm \half} \paran{ \pm \frac{ \varsi{\phi}{i \pm 1} - \varsi{\phi}{i}}{\Delta x}}, \label{Eq:fv-Solid-Current-Flux-nd}
\end{align} 
for ${\rm i} = \mathrm{2,3,\ldots, \params{\mathcal{N}}-1}$. Fluxes at ${\rm i} = 1$ and ${\rm i} = \params{\mathcal{N}}$ are given by the respective boundary conditions. The dimensionless $\varsi{\sigma}{i \pm \half}$ is 1 in active material, $1/\varepsilon_{\sigma}$ in current collector, and a harmonic mean of these values at their interface $\Gamma_{am-cc}$.

\begin{figure}[htbp]
    \centering
    \includegraphics[width=0.9\linewidth]{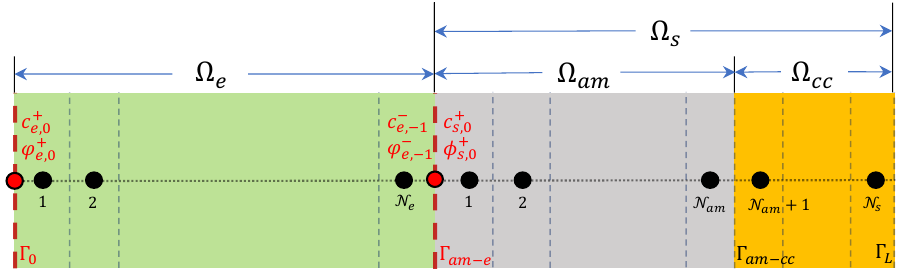}
    \caption{Auxiliary variables introduced at the interfaces $\Gamma_{0}$ and $\Gamma_{am-e}$ for evaluation of the \BV current density.}
    \label{fig:aux_vars}
\end{figure}
  
As detailed above, the discrete molar fluxes and current densities on cell faces are obtained using the cell-centered variables. However, special treatment is needed at interfaces between the solid and electrolyte domain.
The interface conditions at $\Gamma_{0}$ (given in \cref{Eq:BC-Ne-at-interface-0}) and $\Gamma_{am-e}$ (given in \crefrange{Eq:BC-Ne-at-interface-am-e}{Eq:BC-is-at-interface-am-e}) require the \BV current densities which are functions of concentrations and potentials. These variables, however, are stored at the cell centers and need to be evaluated at the interfaces $\Gamma_{0}$ and $\Gamma_{am-e}$. 
A simple way is to extrapolate the adjacent cell-centred variables to the interfaces.
However, this method leads to inaccuracies, especially at higher current rates \cite{zhang_comparison_2014}. 
In this paper, we resolve this problem by employing another approach proposed in \cite{zhang_comparison_2014} involving extra variables at the interfaces.
Thus, we introduce the auxiliary variables, i.e., $\auxvar{c}{e,0}{+}$ and $\auxvar{\varphi}{e,0}{+}$ at $\Gamma_{0}$ and  $\auxvar{c}{e,-1}{-}$, $\auxvar{\varphi}{e,-1}{-}$, $\auxvar{c}{s,0}{+}$ and $\auxvar{\phi}{s,0}{+}$ at $\Gamma_{am-e}$ (see \cref{fig:aux_vars}).
The system of discretized equations is then augmented with an auxiliary system of equations derived from the interface conditions at $\Gamma_{0}$ and $\Gamma_{am-e}$.
This auxiliary system is given by the following relationships between the \BV current density and the gradients of concentrations and potentials
\begin{align}
    G_{e, \rm{c}}^{(A)} &:= \frac{c_{e,1}-\auxvar{c}{e,0}{+}}{\Delta x/2} + \zeta_{e, \rm{c}}\, i_{BV}^{(A)}(\auxvar{c}{e,0}{+} , \auxvar{\varphi}{e,0}{+}) = 0 , \label{Eq:GNp1} \\
    G_{e, \rm{\phi}}^{(A)} &:= \frac{\varphi_{e,1}-\auxvar{\varphi}{e,0}{+}}{\Delta x/2} +  \zeta_{e, \rm{\phi}}(\auxvar{c}{e,0}{+})\, i_{BV}^{(A)}(\auxvar{c}{e,0}{+} , \auxvar{\varphi}{e,0}{+}) = 0 , \label{Eq:GNp2} \\
    G_{e, \rm{c}}^{(C)} &:=  \frac{\auxvar{c}{e,-1}{-}-c_{e,\parame{\mathcal{N}}}}{\Delta x/2} - \zeta_{e, \rm{c}}\, i_{BV}^{(C)}(\auxvar{c}{e,-1}{-}, \auxvar{\varphi}{e,-1}{-}, \auxvar{c}{s,0}{+}, \auxvar{\phi}{s,0}{+}) = 0 , \label{Eq:GNp3} \\
    G_{e, \rm{\phi}}^{(C)} &:= \frac{\auxvar{\varphi}{e,-1}{-}-\varphi_{e,\parame{\mathcal{N}}}}{\Delta x/2} -  \zeta_{e, \rm{\phi}}(\auxvar{c}{e,-1}{-})\, i_{BV}^{(C)}(\auxvar{c}{e,-1}{-}, \auxvar{\varphi}{e,-1}{-}, \auxvar{c}{s,0}{+}, \auxvar{\phi}{s,0}{+}) = 0 , \label{Eq:GNp4} \\
    G_{s, \rm{c}}^{(C)} &:= \frac{c_{s,1}-\auxvar{c}{s,0}{+}}{\Delta x/2} - \zeta_{s, \rm{c}}\, i_{BV}^{(C)}(\auxvar{c}{e,-1}{-}, \auxvar{\phi}{e,-1}{-}, \auxvar{c}{s,0}{+}, \auxvar{\phi}{s,0}{+}) = 0,\label{Eq:GNp5} \\
    G_{s, \rm{\phi}}^{(C)} &:= \frac{\phi_{s,1}-\auxvar{\phi}{s,0}{+}}{\Delta x/2} - \zeta_{s, \rm{\phi}}\, i_{BV}^{(C)}(\auxvar{c}{e,-1}{-}, \auxvar{\varphi}{e,-1}{-}, \auxvar{c}{s,0}{+}, \auxvar{\phi}{s,0}{+}) = 0 . \label{Eq:GNp6}
\end{align}
The $\zeta$'s in \crefrange{Eq:GNp1}{Eq:GNp6} are parameters defined below
\begin{align}
        \zeta_{e, \rm{c}} &= \frac{1-t_+^0}{F\cse{N}}, \quad
        \zeta_{e, \rm{\phi}}(\vare{c}) = \frac{1}{\cse{i}} + \frac{\kappa_D}{F\cse{N}} \frac{1-t_+^0}{\vare{c}},\\
        \zeta_{s, \rm{c}} &= \frac{1}{F\css{N}}, \quad 
        \zeta_{s, \rm{\phi}} = \frac{1}{\css{i}}.
\end{align}
Moreover, the semi-discretized system defined above has only time derivatives of concentrations, while they are absent for the potentials and the auxiliary variables. 
In fact, we can split these variables into differential ($W$) and algebraic ($Z$) variables
\begin{align}
    W &= \curlybrac{\varei{c}{1}, \ldots, \varei{c}{\parame{\mathcal{N}}},  \varsi{c}{1}, \ldots, \varsi{c}{\params{\mathcal{N}}} },\label{Eq:W-definition-aux} \\
    Z &= \curlybrac{\varei{\varphi}{1}, \ldots, \varei{\varphi}{\parame{\mathcal{N}}},  \varsi{\phi}{1}, \ldots, \varsi{\phi}{\params{\mathcal{N}}}, \auxvar{c}{e,0}{+},  \auxvar{\varphi}{e,0}{+}, \auxvar{c}{e,-1}{-}, \auxvar{\varphi}{e,-1}{-}, \auxvar{c}{s,0}{+},  \auxvar{\phi}{s,0}{+} } . \label{Eq:Z-definition-aux}
\end{align}
Consequently, we obtain a system of Differential-Algebraic Equations (DAEs)
\begin{align}
    \Dot{W} &= F(W,Z) , \quad  F: \mathbb{R}^{\mathcal{N}} \times \mathbb{R}^{\mathcal{N}+6} \mapsto \mathbb{R}^{\mathcal{N}} , \label{Eq:DAE-form-D} \\
    0 &= G(W,Z) , \quad  G: \mathbb{R}^{\mathcal{N}} \times \mathbb{R}^{\mathcal{N}+6} \mapsto \mathbb{R}^{\mathcal{N}+6} , \label{Eq:DAE-form-A}
\end{align}
where 
\begin{align}
    F &= \curlybrac{\varei{F}{1}, \ldots, \varei{F}{\parame{\mathcal{N}}},  \varsi{F}{1}, \ldots, \varsi{F}{\params{\mathcal{N}}} }, \label{Eq:F-definition-aux} \\
    G &= \curlybrac{\varei{G}{1}, \ldots, \varei{G}{\parame{\mathcal{N}}},  \varsi{G}{1}, \ldots, \varsi{G}{\params{\mathcal{N}}}, G_{e, \rm{c}}^{(A)},  G_{e, \rm{\phi}}^{(A)}, G_{e, \rm{c}}^{(C)},  G_{e, \rm{\phi}}^{(C)}, G_{s, \rm{c}}^{(C)},  G_{s, \rm{\phi}}^{(C)} }. \label{Eq:G-definition-aux}
\end{align}
The equations corresponding to $W$ and $Z$ lead to the differential equations \cref{Eq:DAE-form-D} and the algebraic constraints \cref{Eq:DAE-form-A}, respectively. 
Further, the algebraic constraints in \cref{Eq:DAE-form-A} are of two types, a discretized elliptic system of equations for the potentials $\varei{\varphi}{i}$ and $\varsi{\phi}{j}$, and a discretized system of nonlinear Robin boundary conditions for the six auxiliary variables.
Moreover, this system of DAEs is in \textit{semi-explicit} or \textit{Hessenberg} form. 
Such systems can also be found in other engineering problems, e.g. low-Mach heterogeneous combustion problems \cite{francois_2024_low_mach}, where interfacial variables may also be used.
In addition, comprehensive numerical investigations illustrated in \cref{sec:index-1-DAE-illustration} show that the system of DAEs corresponding to the LIB problem is of \textit{index}-1~\cite{hairer_solving_1996,petzold_book}.   
Thus, the time integration of the semi-discretized LIB problem is effectively the integration of a system of \textit{index}-1, \textit{semi-explicit} DAEs. Accurate and robust integration of such problems requires specific numerical methods, which are the object of the next section.

\section{Numerical schemes} \label{sec:num_schemes}
In this section, we present various numerical strategies we have used to integrate the stiff system of DAEs. First, we introduce a family of implicit methods suitable for solving DAEs. Then, a multi-domain strategy for the LIB model is presented with explicit and implicit coupling strategies. Finally, an error estimate is formulated to characterize the optimal coupling frequency.

\subsection{Implicit Runge-Kutta methods for DAEs} \label{subsec:implicit-solve-DAE}
Since the 1980s, a number of integration schemes have been designed for DAE problems.
Among others, backward differentiation formulae and stiffly-accurate Runge-Kutta (RK) methods have shown their efficiency and robustness \cite{hairer_solving_1996, ivp_dae_book}.
In our case, the stiffness induced by the fine discretization of diffusion processes and the multiscale physics invites us to contemplate \textit{L-stable} methods.
Based on this requirement, we can select suitable schemes available in the literature.
In the present 1D study, we can afford to use the well-known Radau5 method \cite{hairer_solving_1996}, a 3-\textit{stage} fully-implicit scheme, which is fifth-order accurate for index-1 DAEs and is \textit{L-stable}, having proven especially successful in the integration of stiff ODEs and DAEs.
However, for the larger systems expected in the full-scale 3D LIB simulations, this method can become too costly.   
In that case, other schemes such as the ESDIRK methods \cite{kvaerno_singly_2004} can be used. They, too, exhibit good properties to solve DAEs, while maintaining a lower computational cost \cite{francois_multiphysical_2022,francois_2024_low_mach}.

These advanced methods are in principle high-order generalizations of the implicit Euler method. This first-order method is also \textit{L-stable} and copes well with the \textit{index}-1 DAEs.
Taking advantage of its simplicity, we take it as an illustrative example. Let us assume that the time domain $[0,\ t_{\rm{end}}]$ is divided into $N$ timesteps of equal length $\Delta t=t_{\rm{end}}/N$. Thus, the discrete time values are given by $t_{n+1} = t_{n} + \Delta t$ for $n = \mathrm{0,\ 1,\ \ldots,\ \mathit{N}-1}$, where $t_0=0$ and $t_{N}=t_{\rm{end}}$.
Applying the first-order discretization in time for the differential variables $W$, the DAEs become a system of nonlinear equations given by 
\begin{align} 
    R\paran{W^{n+1}, Z^{n+1}; W^{n}} &= 0, \quad R: \mathbb{R}^{\mathcal{N}} \times \mathbb{R}^{\mathcal{N}+6} \mapsto \mathbb{R}^{2\mathcal{N}+6}, \label{Eqn:nonlinear_problem}
\end{align}
where the residual function $R$, evaluated using previous solution $W^n$, is defined as 
\begin{align}
    R\paran{W^{n+1}, Z^{n+1}; W^{n}} &= 
        \begin{bmatrix} 
            W^{n+1} - W^{n} &-& \Delta t \: F(W^{n+1}, Z^{n+1}) \\
            0 &-& \Delta t \: G(W^{n+1}, Z^{n+1})
        \end{bmatrix}. \label{eqn:residual-funtion}    
\end{align}
Here, the superscript $n$ represents the variables at time $t_{n}$ and so on. The nonlinear problem in \cref{Eqn:nonlinear_problem} can be solved by the Newton-Raphson method, using the solution of the previous timestep (i.e., $W^n$ and $Z^n$) as initial guess. Each Newton iteration requires the solution of a linear system, which can be obtained either by direct or iterative solvers.
Direct solvers are a good choice for 1D and 2D problems, since the bandwidth of the matrices to be inverted remains small. This is not the case in 3D anymore, hence iterative solvers must be considered.
The performance of such solvers is however strongly deteriorated by the ill-conditioning of the matrices obtained, which is partly due to the large number of unknowns in our mixed elliptic-parabolic framework, but also to the multiphysics and multiscale nature of the problem.
It is therefore crucial to use efficient preconditioners (e.g., see \cite{allen_segregated_2021,fang_parallel_2019}) to enable the use of iterative solvers. This topic is beyond the scope of this work and will be treated in a subsequent work devoted to the solution of the full 3D LIB model based on the technique presented in the present paper.

Multirate methods could potentially be envisioned in order to treat the various time scales involved. Whereas they have been mainly developed within the framework of operator coupling such as in \cite{woodward2024} relying on the interesting analysis introduced in \cite{kvaerno2001}, such methods can also be used in coupling various subsystems as in  \cite{constantinescu2023}. However, such an approach, even if it has been largely investigated recently (see \cite{woodward2021} and \cite{constantinescu2013} for example) still lacks flexibility in terms of how the various domains are handled since they are still integrated within the same global time-scheme and make the use of different models in the various domains more difficult to handle. This is always the case, even if some features of the segregated approach can be injected into the multirate technique as in \cite{soszynska2021}.
It is all the more important for our application since the nonlinearity is strong at the coupling interfaces and multi-domain integration will improve the conditioning of the various linear solvers.

We thus now focus on the design of multi-domain techniques at the nonlinear solution level, as an alternative to costly monolithic schemes where the DAEs are integrated simultaneously over the entire computational domain with a single timestep or relying on multirate methods.
Note that such a technique is expected to have a strong impact on the properties of the linearized residual system (formulated among multiple domains) and related linear system solves since it yields an especially improved condition number for the subproblems
that plays a crucial role in the convergence rate of iterative solvers. This statement is investigated in \cref{sec:cnstudy_subsys} for the sake of completeness. 

\subsection{Multi-domain time integration method} \label{subsec:segregated_domain_scheme}
As discussed in \cref{sec:introduction}, the LIB model is a multiphysics as well as multiscale problem, where the \BV current density at solid-electrolyte interfaces introduces a strong nonlinearity. Mathematically, this results in a system of DAEs, which are stiff and costly to solve. The inefficiency and cost of the solution is strongly impacted by the size of the fully coupled problem. Hence, we propose to use a multi-domain technique splitting the monolithic problem into smaller decoupled subproblems.
Such methods have been implemented and studied for the fluid-structure interaction simulations \cite{tezduyar_modelling_2007, takizawa_multiscale_2011}. In the context of LIB simulations, the domain splitting techniques have been used before in \cite{iliev_domain_2016, iliev_splitting_2017, allen_segregated_2021} but remain essentially first-order in time and use a fixed time step. This may hinder the computational performance and accuracy of the simulation tools.
In the present work, we split the LIB model in a manner similar to \cite{allen_segregated_2021}, with the necessary modifications needed for a finite volume method containing the auxiliary system. 
To improve on the previous bottlenecks, we propose to couple the domains in time via a high-order adaptive coupling strategy, known as \textit{multistep coupling}, initially introduced in a different context \cite{francois_multiphysical_2022,francois_multistep_2023} and prototyped in the Python library \textsc{Rhapsopy}~\cite{Rhapsopy2023}. We encourage interested readers to explore the tutorials therein to understand the implementation of the coupling strategy used in this study through simpler test problems.
The originality of the present contribution relies on the extension of this strategy to couple the LIB subproblems efficiently as well as its performance analysis. 
Moreover, we also study the impact of this splitting on the conditioning of the linear systems corresponding to the subproblems (cf. \cref*{sec:cnstudy_subsys}).

We split the fully coupled LIB problem into two subproblems corresponding to the electrolyte and solid domains.
The auxiliary variables and their corresponding equations are also split between the two subproblems. Hence, \cref{Eq:GNp3,Eq:GNp4} corresponding to $\auxvar{c}{e,-1}{-}$ and $\auxvar{\varphi}{e,-1}{-}$ as well as \cref{Eq:GNp5,Eq:GNp6} corresponding to $\auxvar{c}{s,0}{+}$ and $\auxvar{\phi}{s,0}{+}$  become part of the electrolyte and solid subproblems, respectively. 

The electrolyte and solid domains are coupled by the exchange current density at $\Gamma_{am-e}$, given by the \BV model. Hence, the coupling dynamics of the LIB cell simulation is described by the temporal evolution of the auxiliary variables at $\Gamma_{am-e}$.
For a segregated/multi-domain scheme, we refer to these variables as the coupling variables $U_{c}$, i.e., $U_{c} = \bigbrac{\begin{smallmatrix} \vare{U} & \vars{U} \end{smallmatrix} }^T = \bigbrac{\begin{smallmatrix} \auxvar{c}{e,-1}{-} & \auxvar{\varphi}{e,-1}{-} & \auxvar{c}{s,0}{+} & \auxvar{\phi}{s,0}{+} \end{smallmatrix} }^T$. The nonlinear system of  \crefrange{Eq:GNp3}{Eq:GNp6} that corresponds to the coupling variables is referred to as the synchronizing system $G_{\rm{sync}}$.

As for the fully coupled discrete system, each subproblem contains a system of DAEs, which is integrated using the high-order, adaptive implicit methods (described in \cref{subsec:implicit-solve-DAE}). Thus, each subproblem is associated with a residual function whose definition depends on the integration scheme being used. 
In general, the residual functions related to each subproblem can be written as, 
\begin{align}
    \vare{R}(\vare{W}, \vare{Z}, \vars{U}) &= 0 , \quad \vare{R}: \mathbb{R}^{\mathcal{N}_e} \times \mathbb{R}^{\mathcal{N}_e+4} \times \mathbb{R}^{2} \mapsto \mathbb{R}^{2\mathcal{N}_e+4} ,  \label{Eqn:electrolyte_residual} \\ 
    \vars{R}(\vars{W}, \vars{Z}, \vare{U}) &= 0 , \quad \vars{R}: \mathbb{R}^{\mathcal{N}_s} \times \mathbb{R}^{\mathcal{N}_s+2} \times \mathbb{R}^{2} \mapsto \mathbb{R}^{2\mathcal{N}_s+2} , \label{Eqn:solid_residual}
\end{align}
where the subscripts $e$ and $s$ correspond to the electrolyte and solid, respectively.
When \cref{Eqn:electrolyte_residual,Eqn:solid_residual} are solved independently, we need to exchange the required coupling variables, $\vare{U}$ and $\vars{U}$, at regular coupling intervals as done in \cite{allen_segregated_2021}. However, this classical strategy freezes the value of the coupling variables during a coupling interval.
This limits the overall accuracy of the numerical integration to first-order, even if high-order integrators are used in the subproblems. A  way to improve this strategy is to define and use polynomial-in-time approximations to evaluate the coupling variables, so that the integration of each subproblem can be performed with a more accurate evolution of the coupling variables. The polynomials are built by interpolating the values of $U_c$ at the previous coupling time points.
It can be easily shown, following the lines of \cite{francois_multiphysical_2022,francois_multistep_2023}, that for a polynomial approximation of degree $p-1$ and a coupling interval of size $\Delta t_c$, the coupling error $\epsilon_c$ follows
\begin{align}
    \epsilon_c &\sim \mathcal{O}(\Delta t_c^{p}) . \label{Eq:cosim_dtc_law}
\end{align}
We remark that it is important to study the validity of \cref{Eq:cosim_dtc_law} for the 1D LIB half-cell simulations. We have performed this analysis and show our results in \cref{subsec:cosim_order_study}.

The multi-domain technique with coupling of the subproblems at fixed intervals, is very sensitive to the length of this interval (coupling step size $\Delta t_c$). To better understand this, let us consider in more detail the steps we go through during a coupling timestep for time $t_n$ to $t_n+\Delta t_c$. 
At time $t_n$, the approximation polynomials $\hat{U}_c(t)$ of the coupling variables $U_c(t)$, $t \in [t_n, t_n+\Delta t_c]$, are obtained by extrapolation of the $p$ previous values $U_c(t-i\Delta t_c),~i=0\ldots p-1$.
The nonlinear systems in \cref{Eqn:electrolyte_residual,Eqn:solid_residual} are solved to integrate each subproblem independently.
In these subproblems, the interface conditions are obtained by evaluating $\hat{U}_c(t)$ at the required substeps.
At time $t_{n+1}=t_n+\Delta t_c$, a synchronization step is performed, which effectively exchanges information (transfer of coupling variables) between the subproblems. Finally, the approximation polynomials of the coupling variables are updated with their values obtained from the synchronization step. In an \textit{explicit-coupling} scheme (\cref{alg:explicit_cosim}), the solution obtained at this stage is accepted, and we move to the next coupling interval. However, the explicit nature of the extrapolation process is prone to divergence when the time step is too large, and higher-order extrapolations typically have a lower stability limit \cite{francois_multistep_2023}.
\begin{algorithm}[htbp]
    \caption{Algorithm for \textit{explicit-coupling}}
    \begin{algorithmic}
    \STATE{Setup 1D mesh $\mathbf{X}:\mathbb{R}^{\mathcal{N}}$}
    \STATE{Get initial solution states, $W^0$ and $Z^0$}
    \STATE{Initialize approximation polynomial $U_c$}
    \WHILE{$t_n < t_{\rm{end}}$}
    \STATE{$[\vare{W}^{n}, \vars{W}^{n}] \gets W^{n}$ and $[\vare{Z}^{n}, \vars{Z}^{n}] \gets Z^{n}$}
    \STATE{Solve $\Omega_e$: $(\vare{W}^{n+1}, \vare{Z}^{n+1}) \gets R_e(\vare{W}^n, \vare{Z}^n, U_{c}(t))$}
    \STATE{Solve $\Omega_s$: $(\vars{W}^{n+1}, \vars{Z}^{n+1}) \gets R_s(\vars{W}^n, \vars{Z}^n, U_{c}(t))$}
    \STATE{$W^{n+1} \gets [\vare{W}^{n+1}, \vars{W}^{n+1}]$ and $Z^{n+1} \gets [\vare{Z}^{n+1}, \vars{Z}^{n+1}]$}
    \STATE{Synchronization: $U_{c}(t_{n+1}) \gets G_{\rm{sync}}(W^{n+1}, Z^{n+1})$}
    \STATE{Update predictors for $U_{c}(t)$}
    \STATE{$t_{n+1} \gets t_n + \Delta t_c$}
    \STATE{$(W, Z) \gets \texttt{append}(W^{n+1}, Z^{n+1})$}
    \ENDWHILE
    \RETURN $W$, $Z$
    \end{algorithmic}
    \label{alg:explicit_cosim}
\end{algorithm}

This issue can be addressed by employing  small enough coupling intervals such that the stability limit is satisfied. When this is too costly, an alternative is to resort to a modified version of the coupling scheme, which enables an implicit treatment of the coupling variables and hence improves the stability of the coupled simulation. This \textit{implicit-coupling} method, described in \cref{alg:implicit_cosim}, iterates on the definition of $\hat{U}_c$ until the condition $\hat{U}_c(t_{n+1}) \to U_c(t_{n+1})$ is met. We impose a certain tolerance on this limit condition, commonly known as the \textit{waveform-relaxation} tolerance $\rm{WR}_{tol}$.
A fixed-point problem on $U_c(t_{n+1})$ is formulated, that can be solved by simple iterations, dynamically under-relaxed iterations, or Newton iterations. In our case, we only use the simple fixed-point iterations since they perform sufficiently well for the 1D LIB problem.  
\begin{algorithm}[htbp]
    \caption{Algorithm for \textit{implicit-coupling}}
    \begin{algorithmic}
    \STATE{Setup 1D mesh $\mathbf{X}:\mathbb{R}^{\mathcal{N}}$}
    \STATE{Get initial solution states, $W^0$ and $Z^0$}
    \STATE{Initialize approximation polynomial $U_c$}
    \WHILE{$t_n < t_{\rm{end}}$}
    \STATE{$[\vare{W}^{n}, \vars{W}^{n}] \gets W^{n}$ and $[\vare{Z}^{n}, \vars{Z}^{n}] \gets Z^{n}$}
    \STATE{$U_{c}^k \gets U_{c}(t_n)$}
    \WHILE{$\texttt{True}$}
    \STATE{Solve $\Omega_e$: $(\vare{W}^{n+1}, \vare{Z}^{n+1}) \gets R_e(\vare{W}^n, \vare{Z}^n, U_{c}^k(t))$}
    \STATE{Solve $\Omega_s$: $(\vars{W}^{n+1}, \vars{Z}^{n+1}) \gets R_s(\vars{W}^n, \vars{Z}^n, U_{c}^k(t))$}
    \STATE{$W^{n+1} \gets [\vare{W}^{n+1}, \vars{W}^{n+1}]$ and $Z^{n+1} \gets [\vare{Z}^{n+1}, \vars{Z}^{n+1}]$}
    \STATE{Synchronization: $U_{c}^{k+1} \gets G_{\rm{sync}}(W^{n+1}, Z^{n+1})$}
    \IF{$\paran{\frac{||U_{c}^{k+1}-U_{c}^{k}||}{\rm{WR}_{tol}\times||U_{c}^{k}||+\rm{WR}_{tol}/10}<1}$}
    \STATE{$U_{c}(t_{n+1}) \gets U_{c}^{k+1}$}
    \STATE{$break$}
    \ENDIF
    \STATE{Update approximation polynomial for $U_{c}(t)$}
    \ENDWHILE
    \STATE{Update predictors for $U_{c}(t)$}
    \STATE{$t_n \gets t_n + \Delta t_c$}
    \STATE{$(W, Z) \gets \texttt{append}(W, Z)$}
    \ENDWHILE
    \RETURN $W$, $Z$
    \end{algorithmic}
    \label{alg:implicit_cosim}
\end{algorithm}

The \textit{implicit-coupling} method however increases the computational cost, which is proportional to the number of \textit{fixed-point} iterations performed overall. It can be shown that the contraction factor associated with the fixed-point iteration is proportional to $\Delta t_c$, hence a smaller coupling timestep yields a faster convergence. If we dynamically adapt the coupling timestep, we may therefore both improve the accuracy of the simulation and lower the number of iterations per timestep. In case of the explicit coupling, a dynamic coupling timestep also helps to maintain the stability of the coupled computation.
Below, we discuss how the coupling timestep can be dynamically selected.

\subsection{Error estimate and adaptive coupling} \label{subsec:err_estimate_adaptive_coupling}
We propose an error estimation method based on the coupling dynamics of the segregated 1D LIB problem. The error in the coupling variables determines the coupling frequency by guiding the selection of an adaptive coupling step size. In the multi-domain scheme, if the subproblems are integrated with sufficiently small tolerances, the global error is primarily dominated by the coupling error, following the approach in \cite{duarte_splitting_2012}. Consequently, the coupling error, evaluated locally through the coupling variables, effectively governs the overall simulation error.  

\Cref{Eq:cosim_dtc_law} shows that the accuracy of coupling variable approximations in a multi-domain simulation depends on the polynomial degree used. Since the exact error cannot be directly evaluated, we estimate it by comparing the results from two successive degrees of polynomial approximation, using the same coupling step size $\Delta t_{c,\rm{1}}$. Therefore, the difference between the coupling variables $\hat{U}_c^{p-1}$ and $\hat{U}_c^{p}$, obtained from results with polynomial degrees $p-1$ and $p$, provides an error estimate
\begin{align}
    \epsilon_p = \norm{\hat{U}_c^{p-1} - \hat{U}_c^{p}}_2 = \mathcal{O}(\Delta t_{c,1}^{p}) - \mathcal{O}(\Delta t_{c,\rm{1}}^{p+1}) &\approx \alpha \Delta t_{c,\rm{1}}^{p} , \label{Eq:err_p}
\end{align}
where $\norm{~\cdot~}_2$ is the $\ell^2$-norm and $\alpha>$~0 is a proportionality constant.
This error is proportional to $\Delta t_{c,{1}}^{p}$, meaning that increasing the polynomial degree reduces the error. To balance accuracy and computational efficiency, we define an adaptive coupling step size $\Delta t_{c,{opt}}$ based on a user-specified error tolerance $tol$ such that $\epsilon_p \leq tol$. Considering the case when $\epsilon_p \approx tol$ and \cref{Eq:err_p}, we can express
\begin{align}
    tol = \alpha \Delta t_{c,{opt}}^p. \label{Eq:err_opt}
\end{align}
Rearranging this relation and combining it with our initial error estimate in \cref{Eq:err_p}, we obtain the optimal coupling step size
\begin{align}
    \Delta t_{c,{opt}} = \Delta t_{c,\rm{1}} \left(\frac{tol}{\epsilon_p}\right)^{{1}/{p}}. \label{Eq:t_opt}
\end{align}
This strategy dynamically adjusts the coupling step size to maintain accuracy while optimizing computational efficiency. By ensuring that the error remains within a predefined tolerance, the method also helps maintain numerical stability. If an instability grows, the estimated error increases, leading to a time step reduction that restores the multistep coupling scheme to its stability domain. Later, we evaluate the effectiveness of this adaptive strategy for simulating our 1D LIB half-cell model.

\section{Results and discussions} \label{sec:results_discussions}
In this section, numerical studies are conducted based on simulations of the 1D LIB half-cell model presented in \cref{sec:lib_intro}, using the physical parameters listed in~\cref{tab:constants_parameters}. We recall  that all the transport parameters (conductivities, diffusion coefficients and transference number) are assumed to be constant. This simplification does not alter the numerical complexity of the 1D LIB problem; the strong nonlinearity at the solid-electrolyte interface and the ill-conditioned system due to the multiscale physics still remain.
We use a Python code implementation with suitable packages to perform the numerical simulations and analyses. 
Throughout this section, the numerical integration of the system of DAEs is performed using the Radau5 method. The initial conditions for the differential variables are obtained from \cref{Eq:IC-cs}, while the initialization of the algebraic variables follows \cite{vieira2001direct}, using a transient implicit Euler integration for a negligible time $t=\mathrm{1\times10^{-10}s}$. Preliminary studies such as numerical verification of the LIB half-cell model and accuracy of the spatial discretization have been conducted as a natural preliminary step and are presented in \cref{sec:validation_result,sec:space_conv}, respectively. They are important to build our confidence in the 1D LIB model as well as its numerical implementation before proceeding with the assessment of the numerical strategy presented below. Furthermore, to ensure the reproducibility of our implementation, we provide our code along with all necessary parameters, initial data, and other relevant information in \cite{lib1dmdWebpage}. The Jupyter notebooks therein can reproduce results presented in the following.

\subsection{Accuracy of coupling for multi-domain methods} \label{subsec:cosim_order_study}
In this study, we analyze the accuracy of coupling strategies for the multi-domain methods. 
The coupling flux, i.e. the \BV current density, remains constant in time during the CC mode of operation in our 1D setting. Hence, it is not ideal for studying the coupling strategies as the coupling error is negligible for constant coupling conditions. This would not necessarily be the case in higher-dimensional models.
Hence, in contrast with the previous studies, the multi-domain strategy is effectively studied for a constant voltage (CV) operating mode.
Several simulations using the multi-domain method with an increasing number of coupling intervals $\mathcal{N}_t$ (i.e. a decreasing size of the coupling timestep $\Delta t_c$) are performed for both the \textit{explicit} and \textit{implicit} coupling strategies. 
The initial cell voltage $U=$~0.278~V is obtained after 11 s of an initially fully coupled simulation in CC mode at 1C rate. The CV mode of operation using the multi-domain method begins at $t_{\rm{ini}}=$~11~s and is maintained until $t_{\rm{end}}=$~101~s. We use a finite volume grid of size $\mathcal{N}=200$, with $\mathcal{N}_e=\mathcal{N}_s=100$.

To analyze the accuracy of the coupling strategy, we evaluate the coupling error at $t=t_{\rm{end}}$, defined as  
\begin{align}
    {\rm{error}} =  \frac{\norm{y_{\rm{coupling}} - y_{\rm{ref}}}_2 }{\norm{y_{\rm{ref}}}_2} , \label{Eq:cosim_error}
\end{align}
with $\norm{~\cdot~}_2$ being $\ell^2$-norm. Here $y_{\rm{coupling}}$ is the solution obtained using the multi-domain scheme with coupling at fixed intervals. The reference solution $y_{\rm{ref}}$ is obtained on the same finite-volume grid by integrating the coupled system in a monolithic manner with Radau5 and a fine error tolerance of $\mathrm{1\times10^{-12}}$, thus yielding a quasi-exact solution, i.e. with a negligible error from the time integration. 
The multistep nature of the coupling scheme restricts the initial steps to first order ($p=0$), with the order increasing at each step until it reaches the prescribed order of convergence, once enough past points are available to construct a higher-degree approximation polynomial. In a time-adaptive coupled simulation, this is generally not an issue, as in multistep ODE integrators, since the initial steps are usually performed with very small time steps, keeping the accumulated error negligible. However, in this fixed-interval coupling scenario, such limitation can affect the overall order of convergence. To address this, the first four steps (equal to the maximum studied order) are obtained from the reference monolthic solution. With this approach, all methods can directly start at their respective prescribed orders without accumulating first-order errors.

\begin{figure}[htbp]
    \centering
    \includegraphics[width=0.8\linewidth]{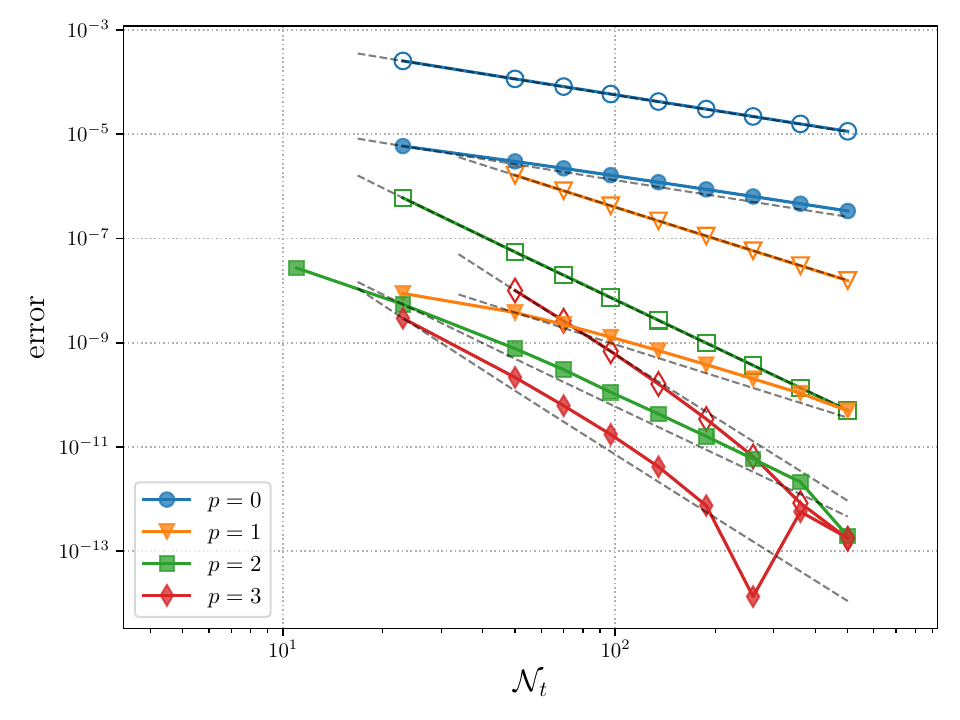}
    \caption{Plots of error vs. $\mathcal{N}_t$ for the multi-domain strategy with coupling at fixed intervals. Errors corresponding to the \textit{explicit} and \textit{implicit} coupling are shown with empty black outlined and filled markers, respectively. The dashed lines represent the theoretical accuracy ($\mathcal{O}(\Delta t_c^{p+1})$) for an approximation polynomial of degree $p$.}
    \label{fig:cosim_convergence_both}
\end{figure}

In \cref{fig:cosim_convergence_both}, the coupling errors from the above simulations are plotted against the number of coupling intervals on a \textit{log-log} scale. The slopes of these curves represent the order of accuracy of the multi-domain method. Each of these curves corresponds to a specific degree of the approximation polynomial and either of the two coupling strategies.
We observe that for a given coupling step size, the \textit{implicit} coupling always has lower errors than the \textit{explicit} coupling strategy. We can explain this by looking at the consistency loops of the \textit{implicit} strategy. Except for the first iteration, interpolation instead of extrapolation is used to approximate the coupling variables. The error constant is smaller for interpolation than that of extrapolation resulting in smaller errors for the \textit{implicit} coupling strategy. 
The simulations are performed for the same set of $\mathcal{N}_t$ for both the \textit{explicit} and \textit{implicit} coupling methods across all coupling orders. In \cref{fig:cosim_convergence_both}, missing data points in a curve indicate that the simulation for that particular case did not converge. Based on the observed trend, we conclude that the \textit{implicit} coupling strategy remains robust for larger coupling intervals (or smaller $\mathcal{N}_t$), unlike its \textit{explicit} counterpart.

Moreover, the error curves satisfy the theoretical order curves given by \cref{Eq:cosim_dtc_law}. Thus, we can safely conclude that the theoretical accuracy of coupling is verified for the 1D LIB problem represented by a system of DAEs. It is worth mentioning that this accuracy had only been verified previously for problems involving coupled subsystems of ODEs in the literature~\cite[Section 8.2]{francois_multiphysical_2022}~\cite{francois_multistep_2023}, not coupled DAEs. Upon this verification, we can confidently use the strategy of \cref{subsec:err_estimate_adaptive_coupling} to evaluate the adaptive coupling step size of the multi-domain method proposed for 1D LIB simulations.

\subsection{Adaptive coupling for oscillating operating condition} \label{subsec:res_adaptive_coupling}
To test the adaptive coupling strategy, we consider a test case involving the application of an oscillating voltage on our 1D LIB model, thus inducing temporal oscillations of the coupling variables.
We take a sine wave signal for the applied voltage with the mean voltage $U_{\rm{mean}}=$~0.103~V $\approx U_0(c_{s,\rm{I}})$. 
It is chosen such that the LIB cell switches between the charging and discharging modes during each oscillation of the voltage signal. This signal consists of three oscillations with an amplitude equal to 5\% of the mean voltage.
The adaptive \textit{implicit} coupling simulations are performed with a tolerance $tol = \mathrm{5\times10^{-6}}$ on the error estimate for the optimal coupling step size. Similar to \cref{subsec:cosim_order_study} each subdomain is discretized with 100 cells for the multi-domain simulations with adaptive coupling.

\begin{figure}[htbp]
  \centering
  \includegraphics[width=0.8\linewidth]{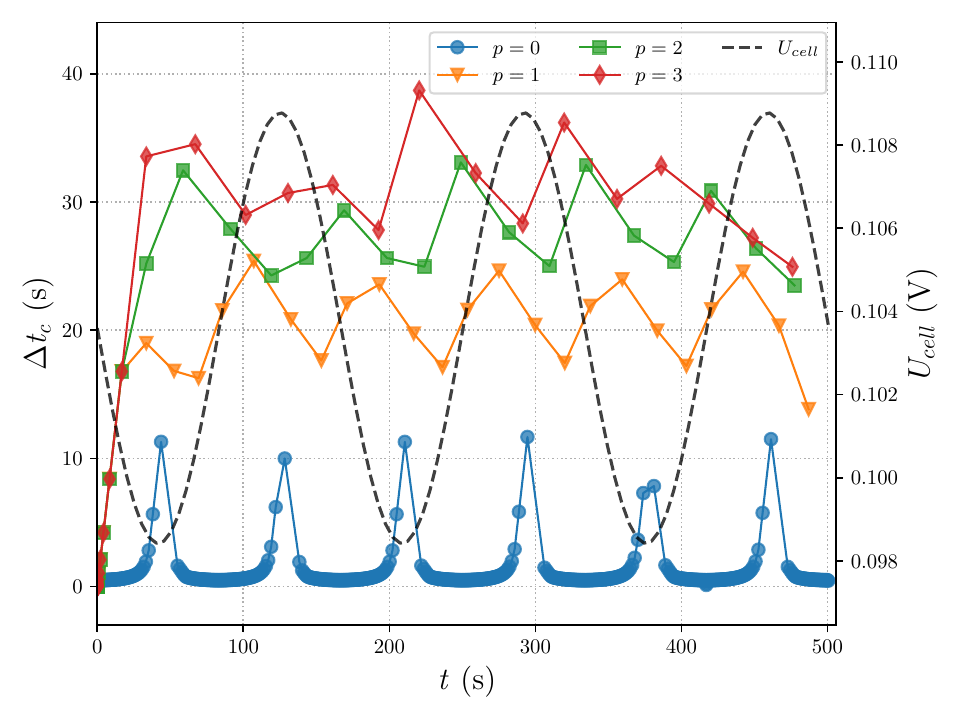}
  \caption{Plots of the coupling timestep $\Delta t_c$~(s) vs. time $t$~(s) for the multi-domain technique with adaptive coupling. Here $p+1$ is the order of the adaptive coupling. The time-varying (sine wave) applied cell voltage is shown with dashed line.}
  \label{fig:sinewave_coupling_dt}
\end{figure}

\begin{figure}[htbp]
    \centering
    \includegraphics[width=0.8\linewidth]{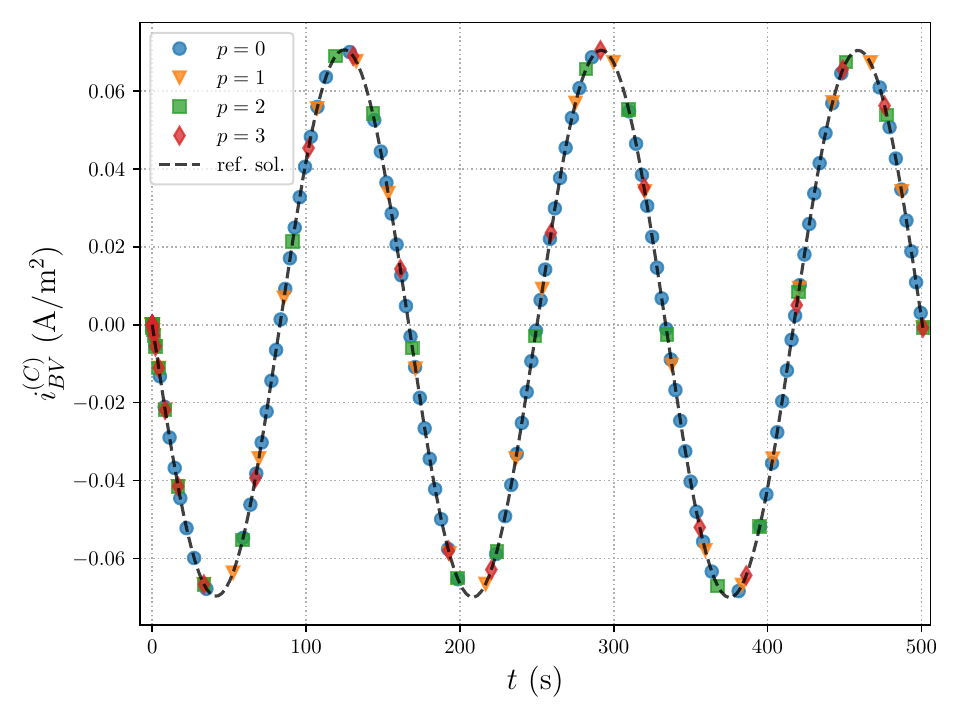}
    \caption{The temporal evolution curves of the \BV current density at $\Gamma_{am-e}$ ($i_{BV}^{(C)}$), evaluated from the simulations using the multi-domain method with adaptive coupling for various coupling orders. For the monolithic reference solution $i_{BV}^{(C)}$ vs. $t$ is shown with dashed line.}
    \label{fig:sinewave_iBV}
\end{figure}

In \cref{fig:sinewave_coupling_dt}, we plot the temporal evolutions of the adaptive coupling timestep $\Delta t_c$ for various orders of the adaptive coupling $p+1$. We observe that $\Delta t_c$ adapts well with the rate of change in the coupling dynamics, which is analogous to the $U_{cell}$ signal.
At the beginning of the simulations, a small coupling timestep value is prescribed. Consequently, as the error reduces, the coupling step size can be increased gradually (with a maximum increase factor of 2 to avoid giving rise to numerical instabilities) until $t\sim$~35~s.
We observe that the simulations with a higher order coupling have larger coupling timesteps,
i.e. for $\mathrm{order} = \curlybrac{\mathrm{1,\ 4}}$ the maximum $\Delta t_c = \curlybrac{\mathrm{11.67~s,\ 38.70~s}}$, respectively.
It may be feared that the larger coupling timesteps can result in a loss of accuracy. 
However, in the present numerical strategy, we ensure that the coupling error, which in turn governs the global error of the simulations, is strictly controlled by the error estimate. Hence, the overall accuracy is maintained even for the larger coupling timesteps. 
We can demonstrate this by comparing the evolution of the \BV current density at $\Gamma_{am-e}$ ($i_{BV}^{(C)}$) in \cref{fig:sinewave_iBV}. For all orders of coupling $p + \mathrm{1}$, there is a good match ($\mathrm{rel.\ error} < tol$) between the simulations with the adaptive coupling and the reference quasi-exact solution. $i_{BV}^{(C)}$ oscillates between positive and negative values indicating the alternating charging and discharging modes of the LIB cell. We also observe that the applied voltage and the interface current density (in \cref{fig:sinewave_coupling_dt,fig:sinewave_iBV}, respectively) are in phase, confirming the absence of capacitive or inductive effects from the mathematical model.

\subsection{Computational performance} \label{subsec:work-precision}
\begin{figure}[htbp]
    \centering
    \includegraphics[width=0.8\linewidth]{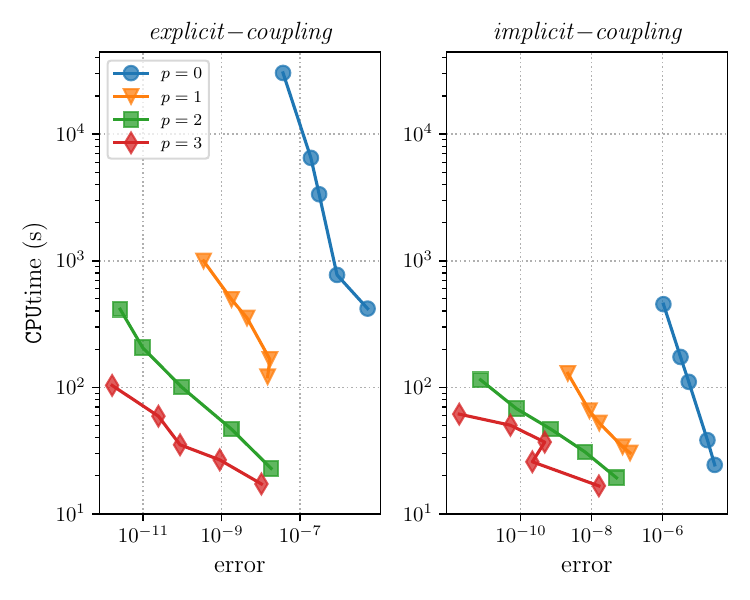}
    \caption{Work-precision diagram for the multi-domain method (\textit{constant voltage operation})}
    \label{fig:constant_WPD}
\end{figure}

\begin{figure}[htbp]
    \centering
    \includegraphics[width=0.8\linewidth]{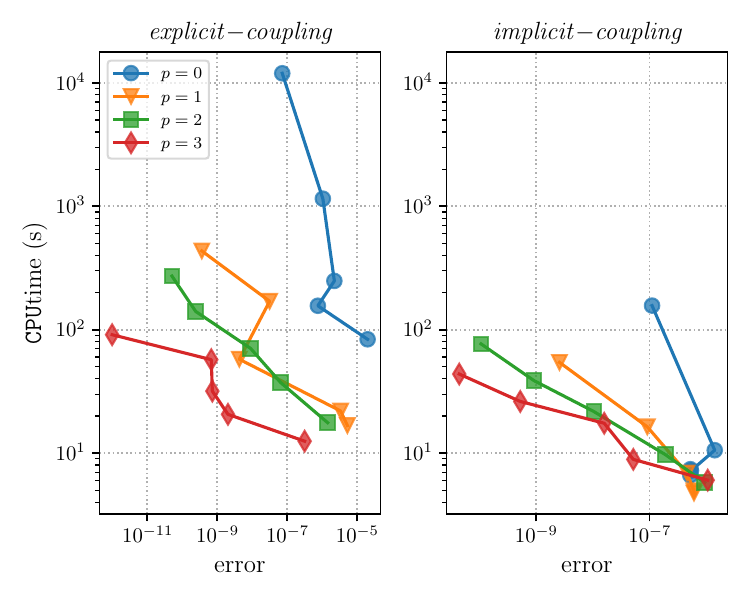}
    \caption{Work-precision diagram for the multi-domain method (\textit{sine wave oscillating voltage operation})}
    \label{fig:sinewave_WPD}
\end{figure}

In this last study, we draw work-precision diagrams to analyze the computational performance of the adaptive coupling strategy implemented in this work. We also tackle the CV case discussed in the accuracy study of \cref{subsec:cosim_order_study}.
These diagrams, shown in \cref{fig:sinewave_WPD,fig:constant_WPD}, illustrate the $\mathtt{CPU} \mathrm{time}$ of various simulations versus their errors evaluated by \cref{Eq:cosim_error} for the different orders of coupling. 
These simulations are performed on a single node of \textit{Cholesky}\footnote{The high-performance scientific computing facility is administrated at Institut Polytechnique de Paris by the IDCS research service unit \url{https://idcs.ip-paris.fr/} - \url{https://docs.idcs.mesocentre.ip-paris.fr/}.} computing cluster with Intel Xeon CPU Gold 6230 20 cores @ 2.1 Ghz.
For each order, the corresponding curve is obtained by performing multiple simulations with increasing values of the coupling error tolerance $tol \in [\mathrm{1\times10^{-10}, 1\times10^{-4}}]$.
The error belonging to each simulation is then computed with respect to a quasi-exact monolithic solution, as done in \cref{subsec:cosim_order_study}.
In most cases, we observe a reduction in $\mathtt{CPU} \mathrm{time}$ as the order of coupling is increased for a given error level. 
Particularly, at higher accuracy levels, a higher-order coupling is even more favorable due to significant performance improvements with respect to lower orders. 
Moreover, the first-order coupling ($p= \mathrm{0}$), which is implemented in previous studies~\cite{iliev_splitting_2017,allen_segregated_2021}, is exceptionally expensive for simulations in both the studied modes of LIB operation. 
We observe that for simulations at $p\geq\mathrm{2}$ both the implicit and explicit coupling perform alike, which might be due to non-stiff coupling dynamics of 1D LIB operation.
The consistency of our observations in this 1D setup suggest the high potential of the proposed numerical strategy for 2D and 3D simulations as well. 

\section{Conclusion} \label{sec:conclusions}
In this study, we have numerically solved a microscale continuum model of LIBs using the finite volume method and an adaptive multi-domain time integration strategy. This choice is motivated by the several orders of magnitude difference in the diffusion time scales characterizing the electrolyte and active material domains.
The main nonlinearity of the problem arises at the solid-electrolyte interface due to the \BV condition modeling the heterogeneous redox reactions. In the discretized system, we introduce auxiliary variables at this interface resulting in a nonlinear Robin condition. Consequently, the subproblem within each domain consists of a stiff system of DAEs
that is solved by the implicit Radau5 method.
The multi-domain approach allows to solve each subproblem independently, while the coupling between the domains is handled by a high-order adaptive coupling scheme.

Our numerical tests have shown that the accuracy of the multi-domain method is governed by the polynomial degree used in the approximation of the coupling variables, and that the adaptive coupling method remains stable even for an oscillating applied voltage at the LIB cell boundaries. Furthermore, higher-order coupling results in larger coupling timesteps, thus reducing the overall computational cost for a desired simulation accuracy.

The results of this 1D study motivate the future implementation of our numerical strategy in more realistic LIB simulations based on 3D electrode microstructures. However, such structures exhibit spatial multiscale features that require a fully adaptive time-space method. Numerical aspects like efficient nonlinear solvers for large systems of DAEs (such as ESDIRK method), efficient preconditioners for the linear solvers, parallelization and related HPC challenges need to be studied in the context of 3D LIB simulations. However, the proposed method can still be envisioned as a breakthrough for these applications in optimizing the coupling timestep for a given accuracy of the solution, while allowing a much better conditioning of the subproblems as illustrated in \cref{sec:cnstudy_subsys}.
Eventually, the multi-domain method offers a suitable framework for the LIB model extension to more complex physics, such as phase transformations and mechanical deformations in the active material, or lithium plating and SEI growth at the solid-electrolyte interface.

\begin{table}[htbp]
    \footnotesize
    \caption{Summary of constants and parameters.}
    \begin{center}
    \ifx\siscFormat\undefined
    \begin{NiceTabular}{|c|c|c|}[hvlines,cell-space-limits=3pt]
    \else
    \begin{tabular}{|c|c|c|}
    \fi
        \hline
        Quantity & Value & Unit \\ 
        \hline
        F & 96487 & C$\cdot$mol$^{-1}$  \\ 
        \hline
        R & 8.314 & J$\cdot$(K$\cdot$mol)$^{-1}$ \\
        \hline
        T & 298.15 & K \\
        \hline
        $L_{e}$ & 20 & $\mu$m \\
        \hline
        $c_{e,\rm{I}}$ & 1000 & mol$\cdot$m$^{-3}$ \\
        \hline
        $D_{e}$ & 1$\times$10$^{-10}$ & m$^2\cdot$s$^{-1}$ \\
        \hline
        $\kappa_{e}$ & 1 & S$\cdot$m$^{-1}$ \\
        \hline
        $t_{+}^0$ & 0.4 & \textit{dimensionless} \\
        \hline
        $\delta_e$ & 0 & \textit{dimensionless} \\
        \hline
        $L_{am}$ & 10 & $\mu$m \\
        \hline
        $c_{s,\rm{I}}$ & 13000 & mol$\cdot$m$^{-3}$ \\
        \hline
        $D_{am}$ & 3$\times$10$^{-14}$ & m$^2 \cdot$s$^{-1}$ \\
        \hline
        $\sigma_{am}$ & 100 & S$\cdot$m$^{-1}$ \\
        \hline
        $F k_0$ & 8.9$\times$10$^{-7}$  & C$\cdot$s$^{-1}\cdot$m$^{2.5}\cdot$mol$^{-1.5}$ \\
        \hline
        $U_0(\vars{c})$ & See Table~2 of~\cite{MAI2020136013} & V \\
        \hline
        $i_{0,\rm{Li}}$ & 10 & C$\cdot$m$^{-2}\cdot$s$^{-1}$\\
        \hline
        $L_{cc}$ & 10 & $\mu$m \\
        \hline
        $\sigma_{cc}$ & 3700 & S$\cdot$m$^{-1}$ \\
        \hline
    \ifx\siscFormat\undefined
    \end{NiceTabular}
    \else
    \end{tabular}
    \fi
    \end{center}
    \label{tab:constants_parameters}
\end{table}

\begin{table}[htbp]
    \footnotesize
    \caption{Summary of characteristic scales.}
    \begin{center}
    \ifx\siscFormat\undefined
    \begin{NiceTabular}{|c|c|c|}[hvlines,cell-space-limits=3pt]
    \else
    \begin{tabular}{|c|c|c|}
    \fi
        \hline
        $\cs{\Phi}$ & $RT/F$ & V\\
        \hline
        $\cs{L}$ & $L_{e} + L_{am} + L_{cc}=L$ & m \\
        \hline
        $\cs{t}$ & $L^{2}/D_{e}$ & s \\
        \hline
        $\cse{c}$ & $c_{e,\rm{I}}$ & mol$\cdot$m$^{-3}$ \\
        \hline
        $\css{c}$ & $c_{s,\rm{max}}$ & mol$\cdot$m$^{-3}$ \\
        \hline
        $\cse{N}$ & $(D_{e}\cse{c})/\cs{L}$ & mol$\cdot$m$^{-2}\cdot$s$^{-1}$\\
        \hline
        $\cse{i}$ & $(\kappa_{e}\cs{\Phi})/\cs{L}$ & C$\cdot$m$^{-2}\cdot$s$^{-1}$\\
        \hline
        $\css{N}$ & $(D_{am}\css{c})/\cs{L}$ & mol$\cdot$m$^{-2}\cdot$s$^{-1}$\\
        \hline
        $\css{i}$ & $(\sigma_{am}\cs{\Phi})/\cs{L}$ & C$\cdot$m$^{-2}\cdot$s$^{-1}$\\
        \hline
    \ifx\siscFormat\undefined
    \end{NiceTabular}
    \else
    \end{tabular}
    \fi
    \end{center}
    \label{tab:scales}
\end{table}


\appendix

\section{Verification of the mathematical model} \label{sec:validation_result}
Here we verify our numerical implementation of the 1D LIB half-cell model by considering a constant current (CC) charge scenario at 0.5C rate. The simulation is performed on a 1D finite volume grid of size $\mathcal{N}=200$. We employ a monolithic time integration method with the Radau5 solver relative tolerance set to $\mathrm{1\times10^{-10}}$. Our results are compared to the reference battery simulator BEST\textit{micro}~\cite{best_software}, using the same physical parameters and grid size. Moreover, numerical results are compared to an analytical solution that we have derived for the CC mode of operation (cf. \cref{sec:analytical_solution}). 

\begin{figure}[htbp]
    \centering
    \includegraphics[width=0.8\linewidth]{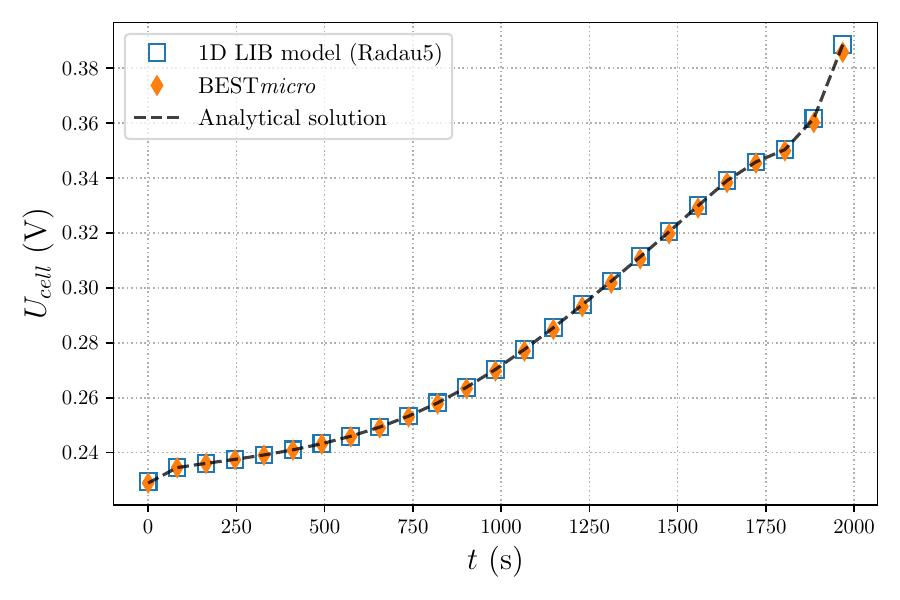}
    \caption{Cell voltage versus time during CC charging at 0.5C rate, obtained from our numerical implementation (current study), BEST\textit{micro}, and analytical solution.}
    \label{fig:validation-1}
\end{figure}

\begin{figure}[htbp]
  \centering
  \includegraphics[width=\linewidth]{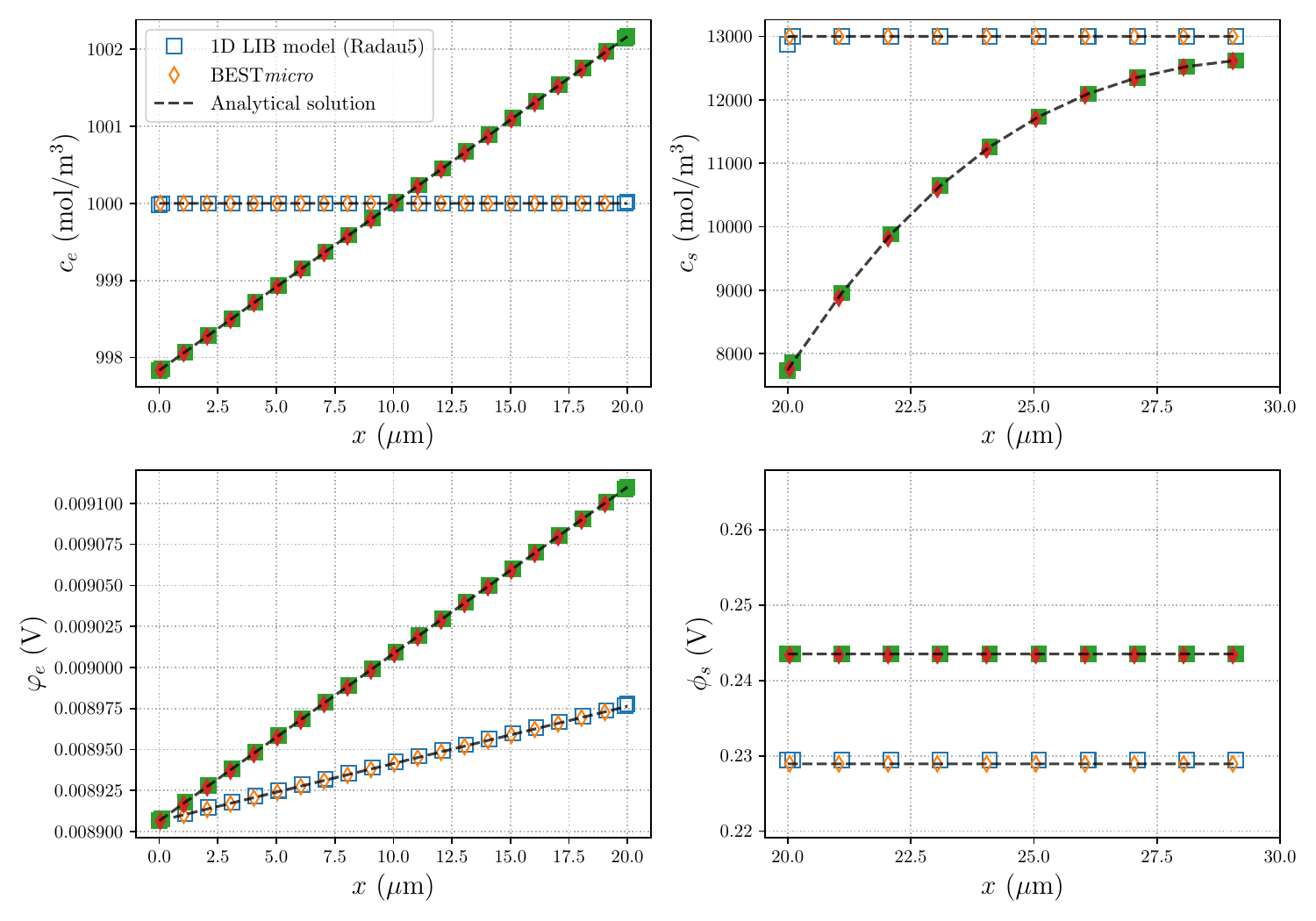}
  \caption{Solution profiles in electrolyte (left panel) and active material (right panel) for CC charging at 0.5C rate, obtained from our numerical implementation (current study) and BEST\textit{micro}. The curves with empty and filled markers represent solutions at $t=$~0~s and at $t=$~500~s, respectively. The dashed lines show the corresponding analytical solution profiles.}
  \label{fig:validation-2}
\end{figure}

In \cref{fig:validation-1}, we plot the temporal evolutions of cell voltage ($U_{cell}$) obtained with our LIB model implementation, BEST\textit{micro}, and the analytical solution. 
Further, in \cref{fig:validation-2}, the spatial solution profiles are displayed at $t=$~0~s and $t=$~500~s for the three cases. These two figures illustrate the good agreement between the solutions, both in time and space. 
Therefore, this benchmark allows us to verify our implementation.       

\section{Space convergence study} \label{sec:space_conv}
\begin{figure}[htbp]
    \centering
    \includegraphics[width=\linewidth]{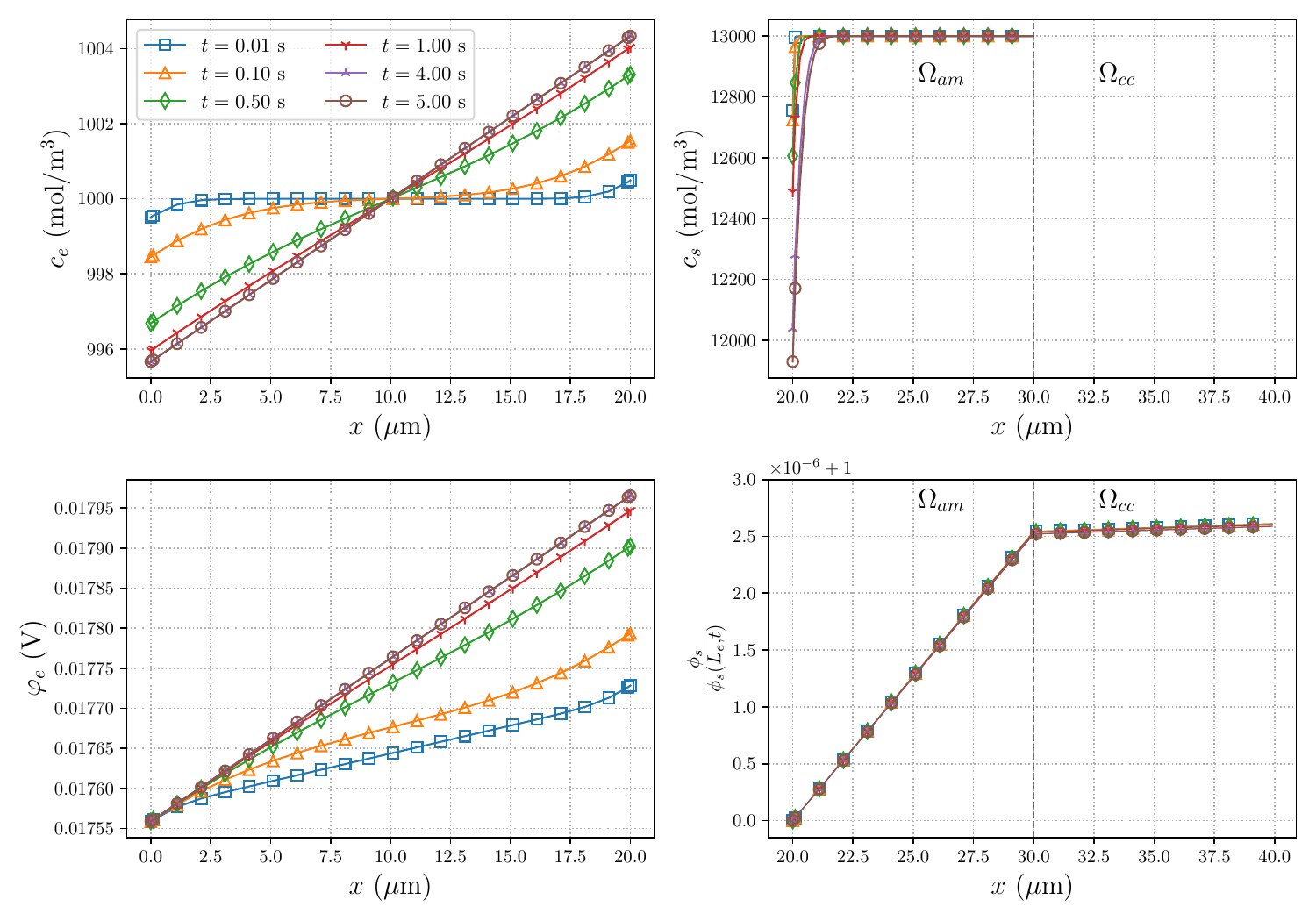}
    \caption{Solution profiles in electrolyte (left panel) and active material (right panel) at various times during CC charging at 1C rate, obtained from our numerical implementation (current study).}
    \label{fig:spatial_profile}
\end{figure}

Using a grid convergence study, we can verify the accuracy of the finite volume discretization. 
We analyze the spatial accuracies of lithium concentrations in the electrolyte and solid domains for a charging scenario at 1C rate. We also study the accuracy of the ionic potential in the electrolyte.
For the electric potential, as we can see in \cref{fig:spatial_profile}, the piecewise linear nature of its solution profiles remains unchanged with time. Only the vertical intercepts (i.e., $\phi_s(L_e,t)$) of these profiles evolve in time due to the variation of the auxiliary electric potential at $\Gamma_{am-e}$. Therefore, we have excluded the electric potential in the solid from this convergence study. 

\begin{figure}[htbp]
    \centering
    \includegraphics[width=0.8\linewidth]{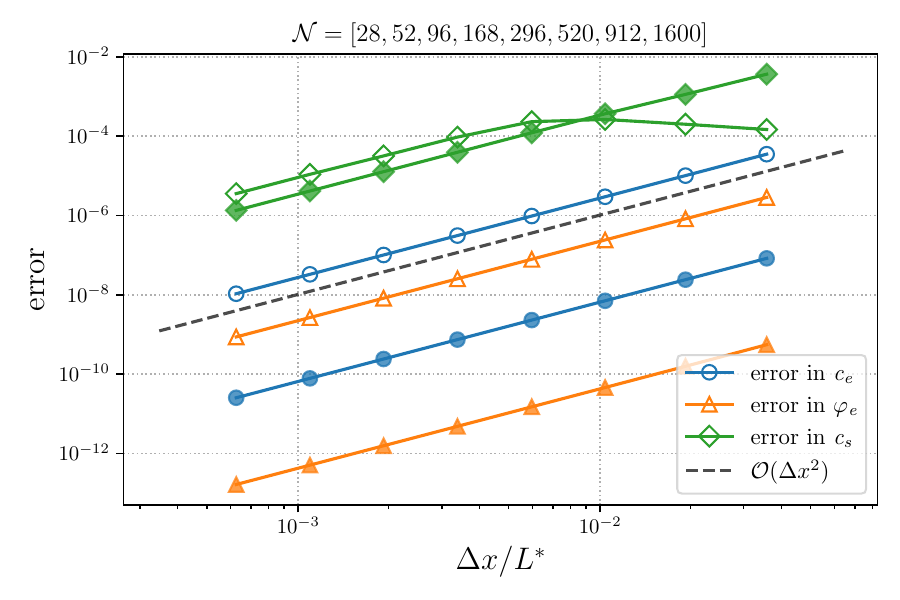}
    \caption{Space convergence plots with errors between the numerical and analytical solutions for charging in CC mode at 1C rate evaluated at $t=$~0.1~s (empty markers) and $t=$~5~s (filled markers).}
    \label{fig:L2-space-conv}
\end{figure}

We perform several simulations with increasing $\mathcal{N}$, i.e., decreasing the 1D cell width $\Delta x$. The spatial integration error is evaluated as
\begin{align}
    {\rm{error}} = \frac{  \norm{y_{\rm{sim}} - y_{\rm{ref}}}_2} {\norm{y_{\rm{ref}}}_2 }  \label{Eq:spatial_error}
\end{align}
where $\norm{~\cdot~}_2$ is $\ell^2$-norm. 
Here $y_{\rm{sim}}$ and $y_{\rm{ref}}$ are the numerical and analytical solution profiles, respectively, both evaluated at the same reference time.
We use sufficiently strict tolerances for time integration using the
Radau5 method such that the temporal discretization errors remain negligible compared to the errors evaluated using~\cref{Eq:spatial_error}.

In \cref{fig:L2-space-conv}, the errors in concentrations and potentials are displayed as a function of dimensionless cell width on a \textit{log-log} scale. The errors are evaluated from solutions at $t=$~0.1~s and $t=$~5~s corresponding to the transient and stationary states of the electrolyte, respectively. We observe that all the curves are parallel to the theoretical second-order curve.
However, the spatial errors in $\vars{c}$ at $t=$~0.1~s are contaminated at larger $\Delta x$, most probably by higher-order terms in the Taylor expansion. These higher-order terms may become significantly large at early times due to the steep spatial gradients of $\vars{c}$ near $\Gamma_{am-e}$.
Nonetheless, we confirm the second order of convergence in space expected from the truncation error of the central difference scheme used in the spatial discretization. 

\section{Analytical solution} \label{sec:analytical_solution}
In the particular scenario of a constant current (CC) charge or discharge, an analytical solution of the 1D LIB half-cell model has been obtained in \cite{romain2022analytical}.
The starting point of the derivation is to observe that the current densities in the electrolyte and solid phases remain constant and are both equal to the imposed current on $\Gamma_{L}$, $i_{\rm ext} = -\xi i_{1C}$. Thus, the electrolyte and solid-phase concentrations are effectively decoupled, each one satisfying a linear diffusion equation with constant Neumann conditions.
The next steps of the derivation are summarized below, using the dimensional notations of \cref{subsec:govn_equations}.

\subsection{Electrolyte concentration and potential} 

Electrolyte concentration satisfies the following initial value problem with Neumann boundary conditions,
\begin{eqnarray} 
    \label{Eq:PDE-varce}
    \frac{\partial\vare{c}}{\partial t}  &=& D_e \frac{\partial^2\vare{c}}{\partial x^2} \ , \\
    \label{Eq:BC-varce}
    \frac{\partial\vare{c}}{\partial x}(0,t) &=& \frac{\partial\vare{c}}{\partial x}(L_e,t) = -\beta_e  \ , \\
    \label{Eq:IC-varce}
    \vare{c}(x,0) &=& c_{e,\rm{I}}  \ ,
\end{eqnarray}
where
\begin{align} \label{Eq:Betae-Definition}
    \beta_e = \left(\frac{1-t_+^0}{F D_e}\right) i_{\rm ext} \ .
\end{align}
Applying the method of separation of variables, we obtain 
\begin{align} \label{Eq:Solution-ce}
    \vare{c}(x,t) &= c_{e,\rm{I}} + \beta_e \left(\frac{L_e}{2}-x\right) \\ \nonumber
    {} & \quad - 4\beta_e L_e \sum_{k=0}^\infty \frac{1}{(2k+1)^2\pi^2} \cos\left(\frac{(2k+1)\pi x}{L_e}\right) e^{-\left(\frac{(2k+1)\pi}{L_e}\right)^2 D_e t} \ . 
\end{align} 
Electrolyte potential satisfies a first-order ODE, obtained by taking the dot product of $\mathbf{i}_e$ from \cref{Eq:Electrolyte-Charge-Conservation} with the 1D normal vector $\mathbf{e}_x$, and then imposing $\vare{\mathbf{i}}\cdot\mathbf{e}_x = i_{\rm ext}$. Further, at $x=0$ the \BV condition \cref{i-BV-A} yields 
$\vare{\varphi}(0,t)=\frac{2RT}{F}\mathrm{argsh}\left(\frac{-i_{\rm ext}}{2 i_{0,\rm Li}}\right)$. By integrating the ODE between $x=0$ and an arbitrary $x$, we obtain
\begin{align} \label{Eq:Solution-phie}
\vare{\varphi}(x,t) = \frac{2RT}{F}\mathrm{argsh}\left(\frac{-i_{\rm ext}}{2 i_{0,\rm Li}}\right) + \frac{2RT(1-t_+^0)}{F}\ln\left(\frac{c_e(x,t)}{c_e(0,t)}\right) - \left(\frac{i_{\rm ext}}{\kappa_e}\right) x \ . 
\end{align} 

\subsection{Solid-phase concentration} 
Let $\vars{\bar{c}}(\xs,t)=\vars{c}(L_e+\xs,t)$ denote the solid-phase concentration in $\Omega_{am}$, expressed in terms of the shifted coordinate $\xs = x - L_e$. This concentration satisfies the following initial value problem with Neumann boundary conditions,
\begin{eqnarray} 
    \label{Eq:PDE-varcs}
    \frac{\partial\vars{\bar{c}}}{\partial t}  &=& D_{am} \frac{\partial^2\vars{\bar{c}}}{\partial \xs^2} \ , \\
    \label{Eq:BC-varcs-1}
    \frac{\partial\vars{\bar{c}}}{\partial \xs}(0,t) &=& -\beta_s  \ , \\
    \label{Eq:BC-varcs-2}
    \frac{\partial\vars{\bar{c}}}{\partial \xs}(L_{am},t) &=& 0 \ , \\
    \label{Eq:IC-varcs}
    \vars{\bar{c}}(\xs,0) &=& c_{s,\rm{I}}  \ ,
\end{eqnarray}
where 
\begin{align} \label{Eq:Betas-Definition}
 \beta_s = \left(\frac{1}{F D_{am}}\right) i_{\rm ext} \ .
\end{align} 
Applying the method of separation of variables, we obtain 
\begin{align} \label{Eq:Solution-cs}
    \vars{\bar{c}}(\xs,t) &= c_{s,\rm{I}} - \beta_s\xs\left(1 - \frac{\xs}{2L_{am}}\right) + \frac{\beta_s D_{am}\,t}{L_{am}} \\ \nonumber
{} & \quad + 2\beta_s L_{am} \left(\frac{1}{6} - \sum_{n=1}^\infty \frac{1}{n^2 \pi^2} \cos\left(\frac{n\pi\xs}{L_{am}}\right) e^{-\left(\frac{n\pi}{L_{am}}\right)^2 D_{am} t} \right) \ . 
\end{align}

\subsection{Solid-phase potential and cell voltage}
Let $c_s^*(t)$ and $\phi_s^*(t)$ denote the solid-phase variables at $\Gamma_{am-e}$, i.e., $c_s^*(t)=\vars{c}(L_e,t)$ and $\phi_s^*(t)=\vars{\phi}(L_e,t)$. Similarly, let $c_e^*(t)=\vare{c}(L_e,t)$ and $\varphi_e^*(t)=\vare{\varphi}(L_e,t)$. Further, we define the cell voltage as $U(t)=\vars{\phi}(L,t)$.
Observing that the solid-phase potential is a continuous and piecewise linear function of $x$ (with constant slopes in $\Omega_{am}$ and $\Omega_{cc}$), we have
\begin{align} \label{Eq:Cell-Voltage-1}
    U(t) &= \phi_s^*(t) - \left(\frac{L_{am}}{\sigma_{am}} + \frac{L_{cc}}{\sigma_{cc}} \right) i_{\rm ext} \ .
\end{align} 
Finally, using the \BV equation \cref{i-BV-C}, we obtain
\begin{align}
    \phi_s^*(t) &= \varphi_e^*(t) + U_{0}(c_s^*(t)) + \frac{2RT}{F}\mathrm{argsh}\left(\frac{-i_{\rm ext}}{2 Fk_0 \sqrt{c_e^*(t)c_s^*(t)\paran{c_{s, \rm{max}}-c_s^*(t)}}}\right) \ , \label{Eq:Cell-Voltage-2}
\end{align}
where $c_e^*(t)$, $\varphi_e^*(t)$ and $c_s^*(t)$ can be evaluated from \cref{Eq:Solution-ce}, \cref{Eq:Solution-phie} and \cref{Eq:Solution-cs}, respectively.

\section{Conditioning of the Jacobian system} \label{sec:cnstudy_subsys}
One of the motivations to use multi-domain methods is to improve the conditioning of the Jacobian system by decoupling the two diffusion time scales. In this study, we evaluate and compare the condition numbers of the Jacobian matrices $J_{\rm full}(W, Z)$, $J_{e}(W_{e}, Z_{e})$, and $J_{s}(W_{s}, Z_{s})$, corresponding to residual functions of the fully-coupled system and to each subproblem (electrolyte and solid), respectively. In this study, we have numerically evaluated these Jacobian matrices for the residual functions obtained from the implementation of Radau5 time integration method.
\begin{figure}[htbp]
    \centering
    \includegraphics[width=\linewidth]{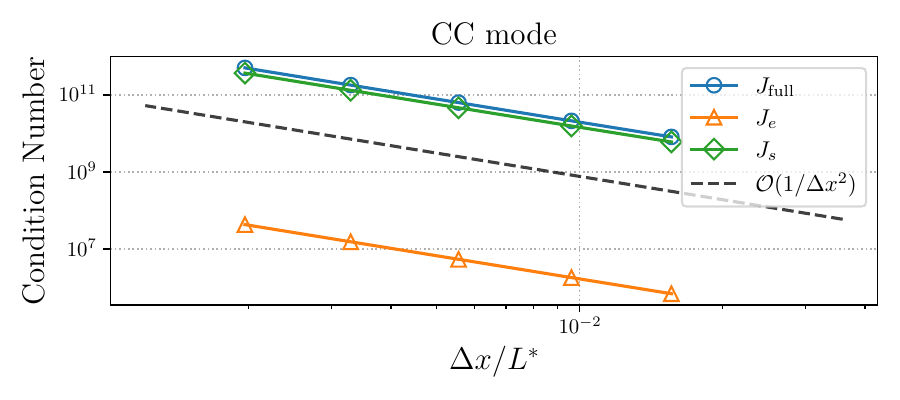}
    \includegraphics[width=\linewidth]{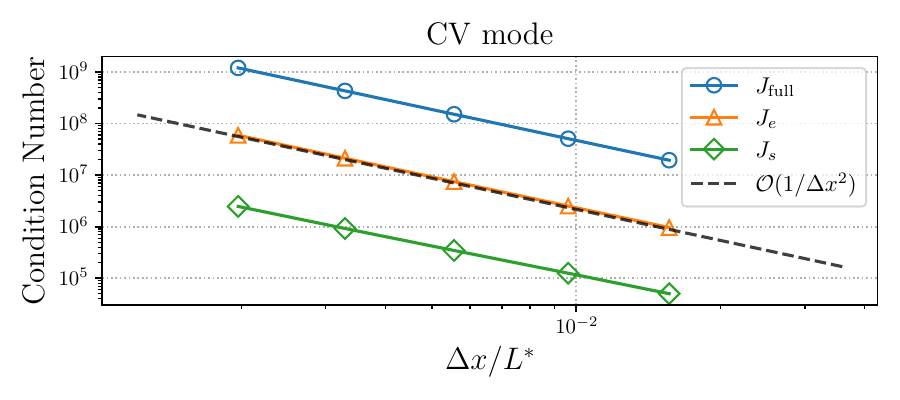}
    \caption{Conditioning of Jacobian matrices $J_{\rm full}$, $J_{e}$, and $J_{s}$, evaluated at $t=\mathrm{101\ s}$ for solutions in CC mode at 1C rate and CV mode at 0.103V.}
    \label{fig:cnStudy_CC_CV}
\end{figure}

\begin{figure}[!htbp]
    \centering
    \includegraphics[width=0.75\linewidth]{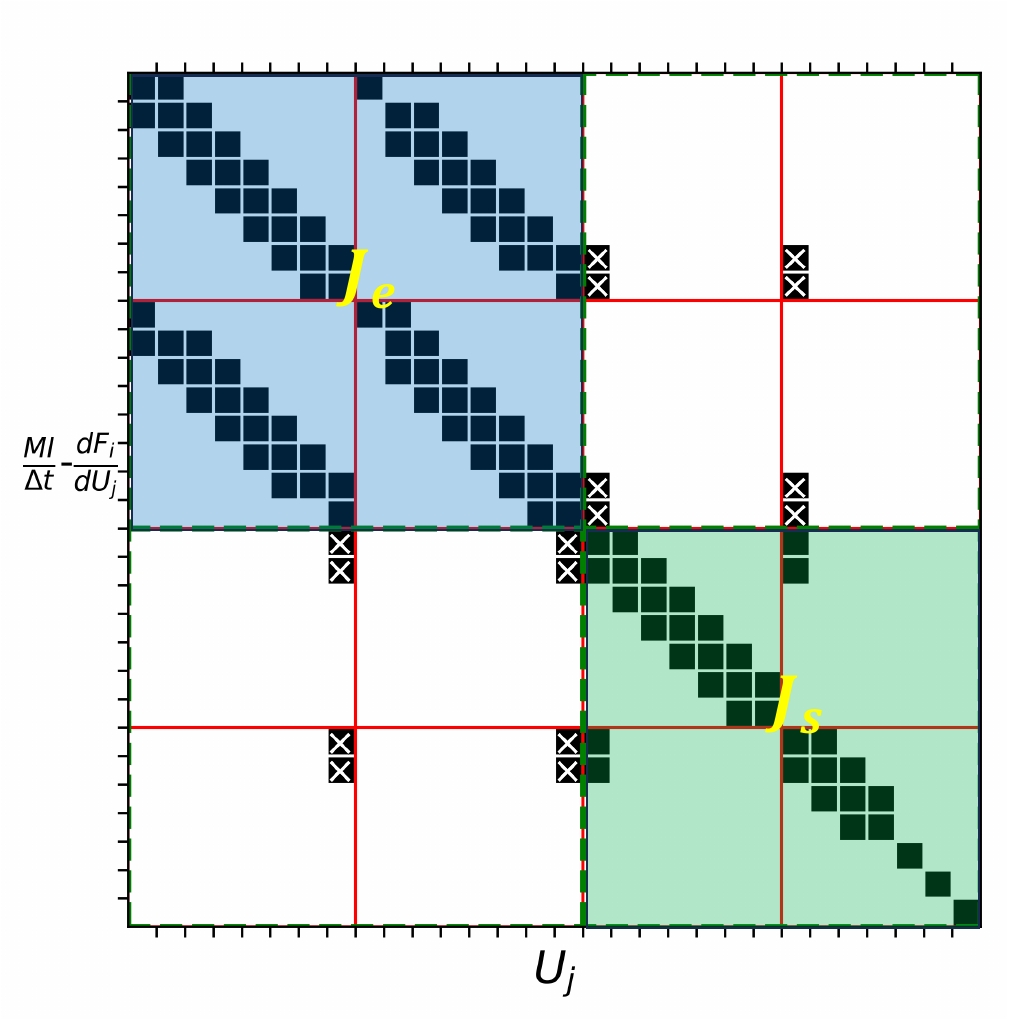}
    \caption{Structure of Jacobian matrices $J_{\rm full}$, $J_{e}$ (shaded blue), and $J_{s}$ (shaded green). The crossed elements are eliminated as a result of decoupling process when using the multi-domain method.}
    \label{fig:residual_jac}
\end{figure}

In \cref{fig:cnStudy_CC_CV}, the condition numbers of $J_{\rm full}$, $J_{e}$, and $J_{s}$ are plotted at the final simulation time $t=\mathrm{101\ s}$ with respect to nondimensional cell width, for CC and CV modes of operation. It should be noted that these condition numbers are almost independent of the temporal evolution of the solution. The plots follow the theoretical order expected from the discretization of the diffusion-based governing equations. The sparsity patterns of these matrices are shown in \cref{fig:residual_jac}, here we observe that $J_e$ and $J_s$, corresponding to the subsystems of the multi-domain method, do not contain the crossed out elements from the monolithic $J_{\rm full}$ as a result of decoupling.    

\begin{figure}[htbp]
    \centering
    \includegraphics[width=\linewidth]{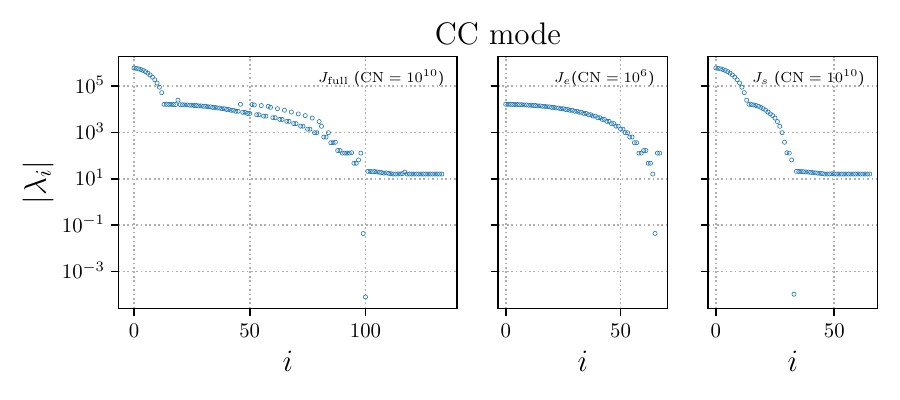}
    \includegraphics[width=\linewidth]{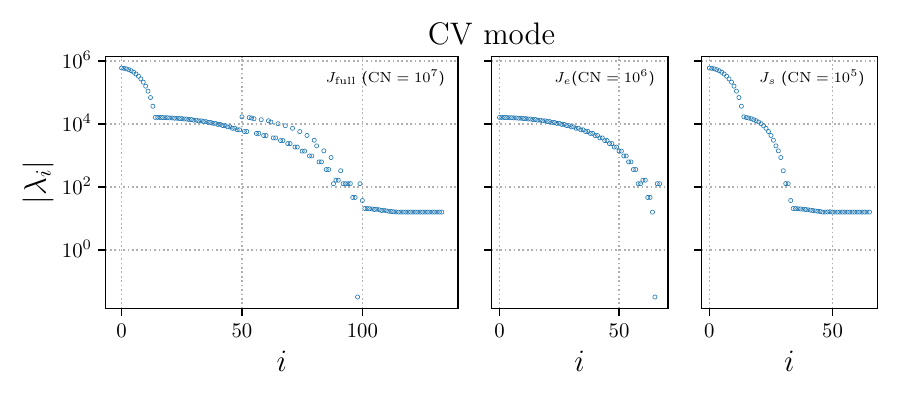}
    \caption{Absolute eigenvalues of the Jacobian matrices $J_{\rm full}$ (left), $J_{e}$ (center), and $J_{s}$ (right) in both CC and CV modes of operation.}
    \label{fig:eig_val_study}
\end{figure}

Further, we investigate the eigenvalues of these matrices in \cref{fig:eig_val_study}. In CC mode, the ill-conditioning of $J_{\rm full}$ corresponds to $J_s$ that has one eigenvalue with magnitude $\ll$1 ($\sim \mathcal{O}(\mathrm{10^{-4}})$). This magnitude is attributed to the Neumann-type boundary conditions on both sides of the solid domain making $J_s$ nearly singular. While in CV mode, we observe a significant reduction in condition number of $J_s$ where the Neumann-type boundary condition is replaced by a Dirichlet-type boundary condition on the right boundary. Such improvements in conditioning are advantageous in the context of iterative linear solvers, to reduce the number of iterations and thus the computational cost. We note that the conditioning of $J_e$ remains unchanged between the two modes of operation as no change is done on this subsystem. 

\section{The DAE description of LIB problem} \label{sec:index-1-DAE-illustration}
We recall from \cref{subsec:semi-discretized} that the semi-discretization of the governing equations results in a system of \textit{index}-1 DAEs. In this study, we show this property via extensive numerical tests. 

\subsection{Index of DAE}
Systems of DAEs are characterized by their \textit{index} $p$, which is the number of times the algebraic constraints, expressed in \textit{semi-explicit} or \textit{Hessenberg} form, need to be differentiated to obtain an explicit expression for $\dot{Z}$. An explicit expression for $\dot{Z}$ will convert the algebraic variables into differential ones, which effectively means converting DAEs into ODEs. Consequently, these ODEs can easily be solved with classical methods while the solution of a system of DAEs of higher index is more complex. Thus, we investigate the index of the system of DAEs corresponding to the microscale LIB model.  

Let us assume that the system of DAEs obtained from the LIB model is of \textit{index}-1. Now we differentiate the algebraic constraints in \cref{Eq:DAE-form-A} with respect to time to get,
\begin{align}
    \paran{\partial_{Z} G}\dot{Z} &= - \paran{\partial_W G}F(W, Z) \ . \label{Eq:Z_dot}
\end{align}
For the \textit{index}-1 assumption to be valid, we need to obtain an explicit expression of $\dot{Z}$. This means that the Jacobian matrix $\partial_{Z} G$ must be invertible. 
We now take a closer look at $\partial_{Z} G$ and show that it is non-singular.

Following from the semi-explicit description of the DAEs in \cref{Eq:DAE-form-D,Eq:DAE-form-A}, we can expand it as follows 
\begin{align}
    \Dot{W} &= F(W, \mathcal{Z}, \mathcal{Z}_{aux}^{(A)}, \mathcal{Z}_{aux}^{(C)}) \ , \ F : \mathbb{R}^{\mathcal{N}} \times \mathbb{R}^{\mathcal{N}} \times \mathbb{R}^{\mathcal{N}_{aux}^{(A)}} \times \mathbb{R}^{\mathcal{N}_{aux}^{(C)}} \mapsto \mathbb{R}^{\mathcal{N}}  \nonumber \\ 
    0 &= \mathcal{G}(W, \mathcal{Z}, \mathcal{Z}_{aux}^{(A)}, \mathcal{Z}_{aux}^{(C)}) \ , \ G : \mathbb{R}^{\mathcal{N}} \times \mathbb{R}^{\mathcal{N}} \times \mathbb{R}^{\mathcal{N}_{aux}^{(A)}} \times \mathbb{R}^{\mathcal{N}_{aux}^{(C)}} \mapsto \mathbb{R}^{\mathcal{N}} \nonumber \\
    0 &= \mathcal{G}_{aux}^{(A)}(W, \mathcal{Z}, \mathcal{Z}_{aux}^{(A)}, \mathcal{Z}_{aux}^{(C)}) \ , \ \mathcal{G}_{aux}^{(A)} : \mathbb{R}^{\mathcal{N}} \times \mathbb{R}^{\mathcal{N}} \times \mathbb{R}^{\mathcal{N}_{aux}^{(A)}} \times \mathbb{R}^{\mathcal{N}_{aux}^{(C)}}\mapsto \mathbb{R}^{\mathcal{N}_{aux}^{(A)}} \ , \nonumber \\
    0 &= \mathcal{G}_{aux}^{(C)}(W, \mathcal{Z}, \mathcal{Z}_{aux}^{(A)}, \mathcal{Z}_{aux}^{(C)}) \ , \ \mathcal{G}_{aux}^{(C)} : \mathbb{R}^{\mathcal{N}} \times \mathbb{R}^{\mathcal{N}} \times \mathbb{R}^{\mathcal{N}_{aux}^{(A)}} \times \mathbb{R}^{\mathcal{N}_{aux}^{(C)}} \mapsto \mathbb{R}^{\mathcal{N}_{aux}^{(C)}} \ ,
    \label{Eq:DAE-form}
\end{align}
where $W$ remains the vector of differential variables while the vector of algebraic variables $Z$ is now expanded into $\mathcal{Z} = \bigbrac{ \begin{smallmatrix} \vare{\varphi} & \vars{\phi} \end{smallmatrix}}^T$, $\mathcal{Z}_{aux}^{(A)}$, and $\mathcal{Z}_{aux}^{(C)}$ such that 
$Z = \bigbrac{\begin{smallmatrix} \mathcal{Z}_{aux}^{(A)} & \vare{\varphi} & \mathcal{Z}_{aux}^{(C)} & \vars{\phi} \end{smallmatrix} }^T $.
Here, $\mathcal{Z}_{aux}^{(A)}$ and $\mathcal{Z}_{aux}^{(C)}$ are the auxiliary variables that were introduced at the solid-electrolyte interfaces of the anode and cathode, respectively. They are defined as,
\begin{align}
    \mathcal{Z}_{aux}^{(A)} &= \{\auxvar{c}{e,0}{+}, \auxvar{\varphi}{e,0}{+} \} \ , \nonumber \\
    \mathcal{Z}_{aux}^{(C)} &= \{\auxvar{c}{e,-1}{-}, \auxvar{\varphi}{e,-1}{-}, \auxvar{c}{s,0}{+},  \auxvar{\phi}{s,0}{+} \} \ .
    \label{Eq:aux-vars-definition}
\end{align}
Consequently, the algebraic constraints are also expanded as 
$G = \bigbrac{\begin{smallmatrix} G_{aux}^{(A)} & \vare{G} & G_{aux}^{(C)} & \vars{G} \end{smallmatrix} }^T $.
Note that for a 1D LIB half-cell, $\mathcal{N}_{aux}^{(A)}=2$ and $\mathcal{N}_{aux}^{(C)}=4$. Also, $\mathcal{N}$ is the size of the discretized problem or number of finite volume cells, with $\mathcal{N}_e$ and  $\mathcal{N}_s$ being the number of cells in the electrolyte and solid domains, respectively.

\begin{figure}[htbp]
    \centering
    \includegraphics[width=\linewidth]{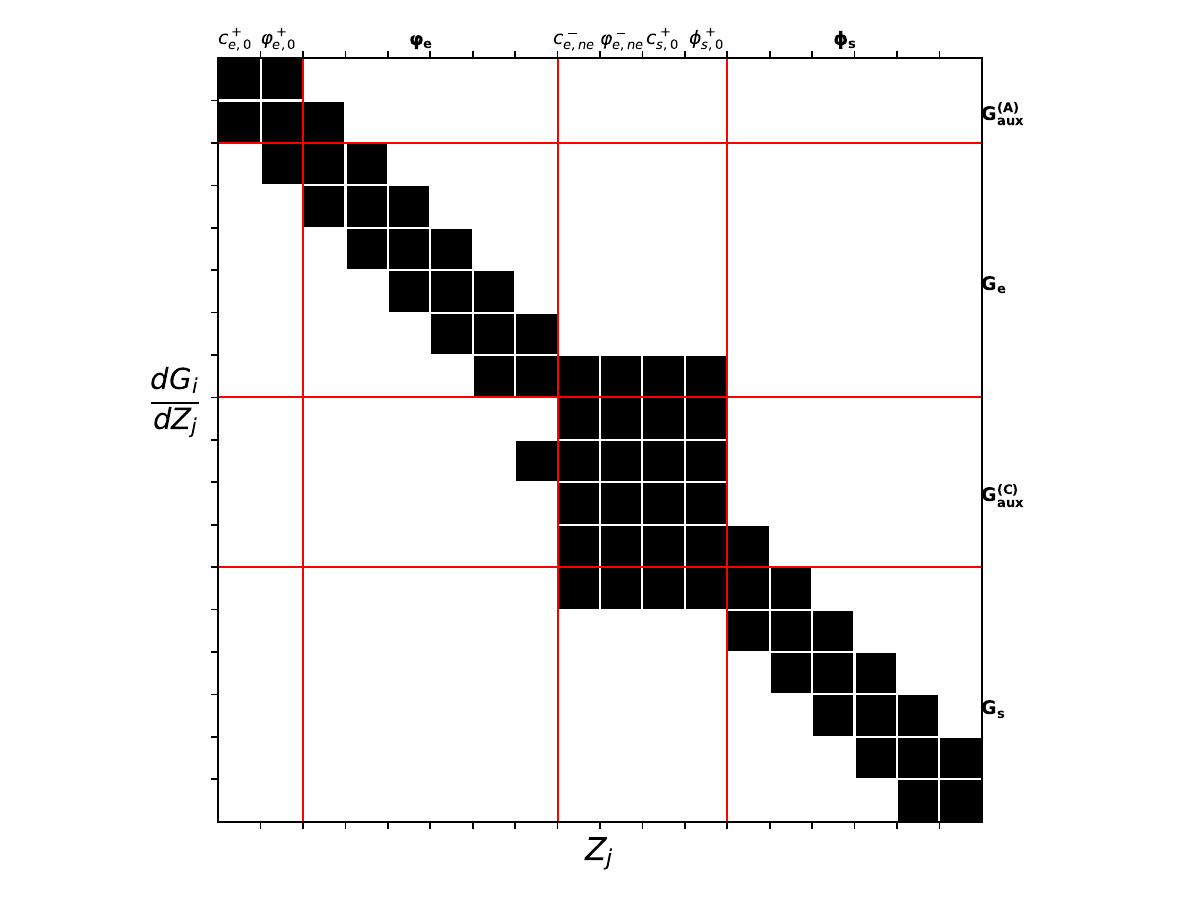}
    \caption{Jacobian of algebraic constraints with respect to algebraic variables $Z$.}
    \label{fig:dae-index-proof}
\end{figure}  
In \Cref{fig:dae-index-proof}, we present the sparsity pattern of the Jacobian $\partial_{Z} G$ or $J_{GZ}$ evaluated on a finite-volume grid with $\mathcal{N}=6$. Further, $J_{GZ}$ can be represented in the block-matrix form,   
\begin{align}
    J_{GZ} &=
    \begin{bmatrix}
        \pder{G_{aux}^{(A)}}{{Z_{aux}^{(A)}}} & \pder{G_{aux}^{(A)}}{\vare{\varphi}} & 0 & 0 \\
        \pder{G_e}{{Z_{aux}^{(A)}}} & \pder{G_e}{\vare{\varphi}} & \pder{G_e}{Z_{aux}^{(C)}} & 0 \\
        0 & \pder{G_{aux}^{(C)}}{\vare{\varphi}} & \pder{G_{aux}^{(C)}}{Z_{aux}^{(C)}} & \pder{G_{aux}^{(C)}}{\vars{\phi}} \\
        0 & 0 & \pder{G_s}{Z_{aux}^{(C)}} & \pder{G_s}{\vars{\phi}}
    \end{bmatrix}, \label{Eq:Jgz_matrix} 
\end{align}
where $J_{GZ}: \mathbb{R}^{\paran{2+\mathcal{N}_e+4+\mathcal{N}_s}} \mapsto \mathbb{M}^{\paran{2+\mathcal{N}_e+4+\mathcal{N}_s}\times\paran{2+\mathcal{N}_e+4+\mathcal{N}_s}}$. We can now use this expression to further investigate the invertibility of $J_{GZ}$.

\subsection{Checking for singularity} \label{subsec:DAE_proofs_nonsigularity}
Firstly, we verify the diagonal dominance property. We note that at rows $3$, $(2 + \mathcal{N}_e)$, $(2 + \mathcal{N}_e + 2)$, $(2 + \mathcal{N}_e + 4)$, and $(2 + \mathcal{N}_e + 4 + 1)$, diagonal dominance is mathematically impossible due the way the physical system is defined. At these rows, there exists one off-diagonal element of same magnitude as the diagonal element.

At rows $3$, $(2 + \mathcal{N}_e)$, and $(2 + \mathcal{N}_e + 4 + 1)$, the diagonal element is equal to $\frac{\mathcal{K}}{\Delta x^2}$ instead of $\frac{2\mathcal{K}}{\Delta x^2}$ due to Neumann-type boundary conditions. The off-diagonal element of the tridiagonal band matrix has the same magnitude. The presence of nonzero elements outside the tridiagonal band structure at these rows (due to the auxiliary system) causes a loss of diagonal dominance. Note that $\mathcal{K}$ is a constant factor that depends on the physical parameters of the problem.

Further, at rows  $(2 + \mathcal{N}_e + 2)$ and $(2 + \mathcal{N}_e + 4)$, the diagonal elements in the diagonal block are $\mathcal{O}(\frac{2}{\Delta x})>>1$. The off-diagonal element outside the diagonal block at these rows has a magnitude equal to $\frac{2}{\Delta x}$. In the working range of our problem, these rows too are never diagonally dominant.    

As we know that the absence of diagonal dominance is not a necessary condition for non-singularity, it is important to investigate other matrix properties such as rank, nature of eigenvalues, condition number etc.

\begin{figure}[htbp]
    \centering
    \includegraphics[width=\linewidth]{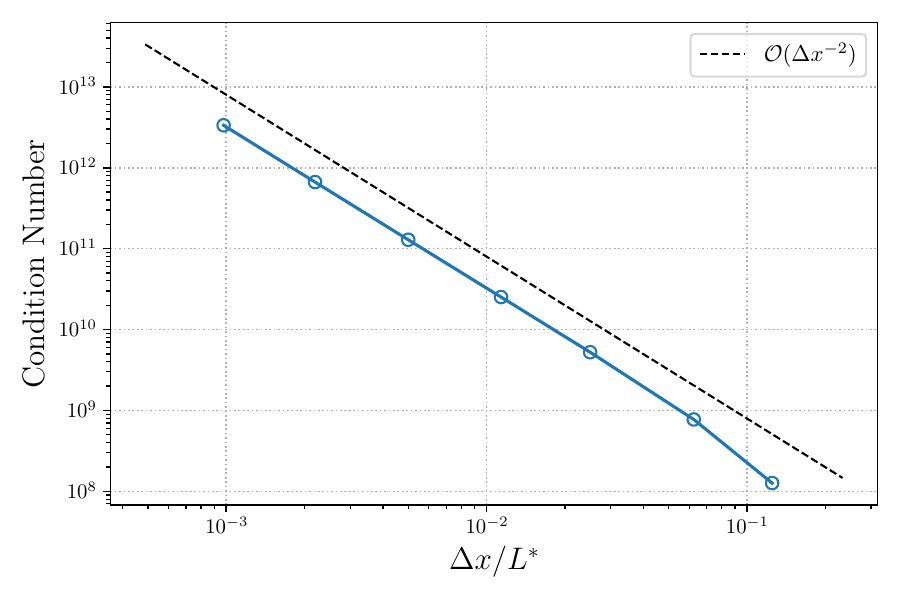}
    \caption{Condition Number of $J_{GZ}$ vs. $\Delta x$.}
    \label{fig:cn-dx}
\end{figure}  
In \Cref{fig:cn-dx}, we observe that the condition number, $CN$ of $J_{GZ}$ is proportional to $\Delta x^{-2}$. We make caution that a large condition number may suggest near singularity of the matrix. However, ill-conditioning is not a necessary condition for the singularity. 

\begin{figure}[htbp]
    \centering
    \includegraphics[width=\linewidth]{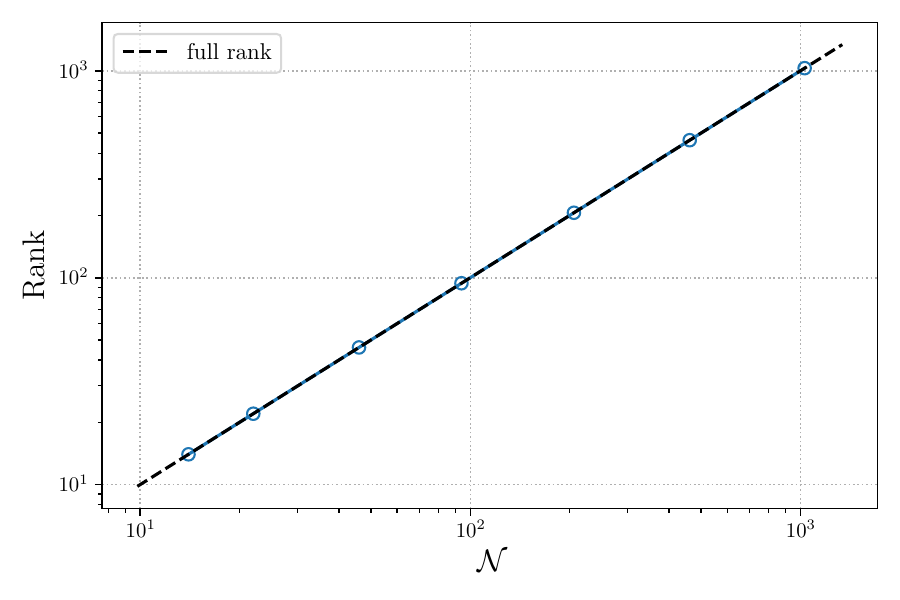}
    \caption{Rank of $J_{GZ}$ vs. $\mathcal{N}$.}
    \label{fig:rank-n}
\end{figure}  
Finally, for each size of system $\mathcal{N}$, the rank of the matrix $J_{GZ}$ is equal to $\mathcal{N}$ (see Figure~\ref{fig:rank-n}), i.e., the Jacobian matrix always has a full rank. We are aware that full rank is a necessary condition for non-singularity and $J_{GZ}$ satisfies it for the tested range of $\mathcal{N}$ and solution vectors $W$ and $Z$. Therefore, considering this comprehensive numerical investigation we can conclude that $J_{GZ}$ is non-singular matrix and hence, \cref{Eq:DAE-form} represent an \textit{index}-1 system of DAEs.

\section*{Acknowledgments}
The authors would like to thank TotalEnergies OneTech for its financial and technical support. They are grateful to their SAFT colleagues in Bordeaux, and their TotalEnergies colleagues at the Playground Paris-Saclay, for useful discussions on battery modeling and simulations. They also acknowledge the valuable help on code development from Laurent Séries and Loïc Gouarin (research engineers at CMAP, École Polytechnique).

\bibliographystyle{plain}
\bibliography{references.bib}


\end{document}